\documentclass[11pt,hidelinks]{article}
\usepackage[english]{babel}
\usepackage{lineno} %line numbers
\usepackage{acronym}
\usepackage{adjustbox}
\usepackage{amsmath}
\usepackage{amssymb}
\usepackage{booktabs}
\usepackage{caption}
\usepackage{csquotes}
\usepackage{courier}
\usepackage[shortlabels]{enumitem}
\usepackage[gen]{eurosym}
\usepackage{floatrow}
\usepackage[letterpaper,top=2.5cm,bottom=2.5cm,left=2.5cm,right=2.5cm]{geometry}
\usepackage{graphicx}
\usepackage{hyperref}
\usepackage{lmodern}
\usepackage{longtable}
\usepackage{makeidx}
\usepackage{multirow}
\usepackage{rotating}
\usepackage{setspace}
\usepackage{siunitx}
\usepackage{subfigure}
\usepackage{titlesec}
\usepackage{titling} % Customizing the title section
\usepackage{wrapfig}
\setlist[itemize]{noitemsep} % Make itemize lists more compact
\setlist[enumerate]{noitemsep}
\title{Mixed Integer Linear Program model for optimized scheduling of a vanadium redox flow battery with variable efficiencies, capacity fade, and electrolyte maintenance}
\author{Cremoncini Diana$^a$*, Frate Guido Francesco$^a$, Bischi Aldo$^a$, Ferrari Lorenzo$^a$}
\date{
    \small{\textit{$^a$University of Pisa, Department of Energy, Systems, Territory and Construction Engineering, Largo Lucio Lazzarino 1, 56122, Pisa, Italy. *Corresponding author\\}}%
    \small{Email addresses: diana.cremoncini@phd.unipi.it (Diana Cremoncini), guido.frate@unipi.it (Guido F. Frate), aldo.bischi@unipi.it (Aldo Bischi), lorenzo.ferrari@unipi.it (Lorenzo Ferrari)}
}

\begin{document}
\maketitle

\rule{15cm}{0.5pt}

\begin{abstract}
    % Background: 
    Redox Flow Batteries are a versatile and durable type of electrochemical storage and a promising option for large-scale stationary energy storage. The vanadium redox flow battery represents the most mature chemistry for the technology, and it is the most widely commercialized system thanks to its high chemical stability and performance.
    %
    % Purpose/aim:
    This work aims to optimize the scheduling of a vanadium redox flow battery that stores energy produced by a renewable power plant, keeping into account a thorough characterization of the battery performance, with variable efficiencies and capacity fade effects.
    A detailed characterization of the battery performance improves the calculation of the optimal number of cycles and revenue associated with the battery use if compared to the results obtained using simpler models, which take into account constant efficiencies and no capacity fade effects.
    %
    % Method:
    The presented problem is nonlinear due to the functions of the battery efficiency, which depend upon charging and discharging powers and state of charge with nonlinear, non-convex correlations.
    Therefore, as per the state of the art of operation research, the problem is linearized using convex hulls in three dimensions.
    The optimization program also calculates the progressive battery capacity fade due to undesired secondary electrochemical reactions and the economic impact of capacity restoration through periodic maintenance.
    The final problem is solved as a Mixed-Integer Linear Program (MILP) to guarantee the global optimality of the linearized problem.
    The proposed optimization model has been applied to two different case studies. The first is a case of energy arbitrage, where the electric energy produced by a renewable plant is stored when the energy price is low and sold to the grid when the energy price is high to maximize the profit. The second is a case of load-shifting, balancing electric energy generation and demand from a grid-connected renewable energy community, where storage minimizes the expenses for the energy purchased from the grid.
    %
    % Results:
    The optimization results have been compared to those obtained with constant battery efficiency models, which do not consider the capacity fade effects. Results show that simpler models overestimate the optimal number of cycles of the battery and the revenue by up to 15\% if they do not take into account the degradation model of the battery, and respectively up to 32\% and 42\% if they also assume constant efficiency for the battery.
\end{abstract}

\rule{15cm}{0.5pt}

\section{Introduction}

    United Nations' Intergovernmental Panel on Climate Change \cite{IPCC2018} explained that the world would have to progressively reduce its carbon dioxide emissions, to limit the global temperature increase to 1.5$^\circ$ C.
    For this reason, the European Union aims to be climate-neutral - having an economy with net-zero greenhouse gas emissions - by 2050, as expressed in the European Green Deal \cite{GreenDeal} and as certified by the EU’s commitment to global climate action under the Paris Agreement \cite{ParisAgreement}.
    To achieve these objectives, it is fundamental to increase the share of energy coming from renewable sources \cite{EU_dir_2018}, and especially for the electricity sector, the rapid expansion of renewable sources integration calls for a more flexible energy system \cite{REN21_2021}. The storage of electrical energy is a prominent solution to address several technical and economic challenges of renewable sources integration and achieve grid reliability and stability \cite{IRENA_stor}.
    
    Redox Flow Batteries (RFBs), a versatile and durable type of electrochemical storage, are emerging as a promising technology for stationary applications and renewable energy support because they offer valuable operational advantages, such as working at ambient temperatures, having a long cyclic life, and having independently sizable power and capacity \cite{IRENA_cost}. This technology could offer competitive cost ranges if used for cases requiring long storage times \cite{Diez2020}.
    
    RFBs use redox-active materials dissolved in liquid electrolytes stored in tanks and pumped through two hydraulic circuits into an electrochemical cell \cite{Sandia_Handbook2015}.
    A semipermeable membrane separates the two sides of the reaction cell.
    An RFB's electrochemical reactor is composed of multiple cells stacked together to increase the overall voltage.
    The battery power depends on cell and stack size, while the battery capacity only depends upon the volume of the electrolytes, which can be increased by simply resizing the storage tanks.
    This modular structure allows for independently sizing energy and power, reaching energy-to-power ratio values of more than 10 hours \cite{Sandia_Handbook2015}.
    RFBs have a long lifecycle, estimated to be greater than 13000 cycles \cite{Alotto2012}, which means ten to twenty years of stable battery operation.
    
    \subsection{Vanadium Redox Flow Battery} \label{VRFB}
    Among RFBs, the all-vanadium (VRFB) is the most researched and technologically successful battery thanks to its high chemical stability and overall performance.
    The VRFB is composed of aqueous acidic solutions of vanadium ions in different oxidation states as active materials in the electrolytes.
    The capital cost of a VRFB is relatively high if compared to other storage technologies: PNNL 2020's report \cite{PNNL_2020} states that for a 100 MW, 10-hour discharge system,
    the mean projected cost for a VRFB is 399 \$/kWh, while the cost of a lithium-ion LFP is 356 \$/kWh.
    
    The standard electrochemical reactions for VRFB charge and discharge are shown for both the negative (on the left) and positive (on the right) half cells in Eq. (\ref{VRFBequation}):
    
    \begin{equation} \label{VRFBequation}
        V^{3+} + e^- \rightleftharpoons V^{2+}
        \quad \hspace{1 cm}
        VO^{2+} + H_2O \rightleftharpoons VO_2^+ + e^- + 2H^+
    \end{equation}
    
    The operational temperature range for this technology is 10-40 $^\circ$C \cite{Bryans2018}. Coupled with nominal cell voltages of 1.1 to 1.6 V, the energy density of the battery stands at about 25-35 Wh/L (noticeably lower than the 250 Wh/L of Li-ion batteries \cite{Diez2020}).
    
    The single cell of a VRFB can reach charging efficiency values over 90\% \cite{Bindner2010}.
    The overall system efficiency is lower due to other losses in the system. Other than intrinsic cell losses \cite{Pugach2019}, there are other losses, due to parasitic shunt current losses in the cell stack, circulation pump, inverter, and the other power control devices; all these losses are influenced by many operational factors, such as current, ionic concentration, temperature, and electrolyte flow rate \cite{Zugschwert2021}.
    The system maximum roundtrip AC/AC efficiency, including all the mentioned stack losses and auxiliary losses, varies from about 60\% \cite{Bindner2010, Zugschwert2021} up to 70\% \cite{Guarnieri2020}, depending on stack and battery management system design.
    
    VRFBs have a long operational life because the reactions involved in the cell do not consume ion metals, avoiding permanently damage of the electrolyte.
    Nevertheless, the battery suffers from degradation effects, which can cause the capacity to fade over time and needs periodic maintenance action to restore the battery capacity to its original value.
    The most important degradation effects occurring in the vanadium battery and relative restoration strategies are the following \cite{Rodby2020,Poli2020,Yuan2019}:
    \begin{itemize}
        \item crossover, an undesirable transport of active material and water through the cell membrane, causes a concentration imbalance, reducing the accessible capacity; this effect is reversed with a periodic electrolyte mixing, followed by a battery recharging, called rebalancing;
        \item oxidative imbalance, caused by undesired oxidation of active species, which also reduces the accessible capacity; this effect is reversed with a periodic chemical or electrochemical maintenance process, called servicing.
    \end{itemize}
    
    To reverse the capacity fade of the electrolytes, mainly caused by two reversible degradation mechanisms, two different maintenance actions are needed. Rebalancing is simple, economical, and can be applied to VRFBs of all sizes paying a price in terms of time when the battery cannot be utilized (in the order of the discharge time of the battery) and electric energy (in the order of 150\% of the rated capacity).
    On the other hand, the servicing process needs dedicated equipment and instruments to be carried out, and its cost cannot be easily identified \cite{Poli2020}.
    \subsection{State of the art}
    Literature suggests many ways to model VRFB systems with different complexity levels to conduct techno-economic studies.
    A battery model can take into account constant or variable system efficiency. The latter could be described via efficiency curves depending on various factors, such as charging/discharging powers, state of charge, temperature, and electrolyte flow rate. It can also take into account various degradation and ageing mechanisms.
    
    As far as the VRFB performance characterization is concerned, Bindner et al. \cite{Bindner2010} point out how the overall system efficiency is far from being constant, while other works \cite{He2016, Jafari2019} show that literature models considering constant battery efficiency in techno-economic assessments \cite{Locatelli2015, Lajimi2019, Terlow2020} lead to relevant errors on the optimal battery management strategy prediction and the economic revenue associated with the battery installation.
    Zhang X. and Skyllas-Kazacos M.\cite{Zhang2016}, define an optimal sizing method of a vanadium redox flow battery for residential use, where the efficiency is variable with charging/discharging power, but they neglect the efficiency dependence on the state of charge and do not take into account the progressive capacity fade.
    Gong et al. \cite{Gong2019}, similarly, define the optimal configuration and scheduling of the energy storage system using efficiencies functions depending on power and state of charge to describe the battery performance, without introducing capacity fade of the electrolyte or other degradation mechanisms.
    Similar optimization models were produced by \cite{Jafari2019} and \cite{Nguyen2019}.
    
    In particular, the work by Jafari et al. \cite{Jafari2019}, consists of creating a detailed Mixed-Integer Linear Program (MILP) for the optimal management of a vanadium flow battery, considering only the internal energy losses of a laboratory cell, without taking into account the auxiliary losses of the battery, such as shunt currents, pumping or inverter losses.
    Jafari et al. include the piece-wise linearization of battery efficiencies, see details in \cite{Jafari2020}, which are expressed as a function of two variables (battery charge/discharge power and state of charge - $SoC$).
    It is assessed in \cite{Jafari2020} that the linearization algorithm can only be applied to convex functions, increasing monotone with respect to one of the two variables.
    Jafari's linearization method, suitable for nonlinear efficiency functions of two variables (operating power and $SoC$), needs the introduction of multiple binary variables to select the state of charge curve, which causes an undesired computational burden to the problem.
    
    The battery model can be further developed and complicated by adding degradation mechanisms, which cause cycle-dependent capacity fade, and, eventually, with the introduction of maintenance methods for capacity restoration.
    The capacity fade model for a VRFB, introduced in a more recent study by Jafari et al. \cite{Jafari2021}, focuses on optimizing the total number of rebalancing events during a whole year and gives the optimal number of maintenance events for a given number of cycles during the battery lifetime.
    Nevertheless, in \cite{Jafari2021}, there is no distinction in the sources of capacity fade, nor a detailed technical description of the maintenance events and their effect on battery utilization.
    As far as other capacity fade models applied to VRFBs, He et al. \cite{He2016} take into account the capacity fade, without introducing any electrolyte maintenance strategy, and solving the problem as non-linear.
    
    \subsubsection{Description of research gaps}
    
    In the framework of studies on VRFBs' optimal scheduling, we want to address and cover the following research gaps:
    \begin{itemize}
        \item battery efficiency is often described as a constant, and when modeled as a variable of other parameters, models are often inaccurate, lacking either system auxiliary losses, or the dependency on the state of charge;
        \item where the efficiency is modeled as a function of power and SoC, the problem is either solved as non-linear, or linearized with the help of multiple binary variables, leading in both cases to a high computational burden;
        \item capacity fade effects or other degradation mechanisms are hardly introduced in optimization models, and when considered they usually contain little details on the technical aspects of the degradation mechanisms, or do not include and historic battery cycling;
        \item detailed electrolyte maintenance processes to restore lost capacity and their impact on optimal scheduling have been considered in detail only in higher level techno-economic studies \cite{Rodby2020}, but to the authors' best knowledge not in terms of technical impact on battery utilization.
    \end{itemize}

    Note that the impact of other parameters on battery efficiency, such as temperature and electrolyte flow rate, have been usually discarded in literature both due to their minor impact on battery performance and the lack of experimental data.
    
    These literature gaps show oversimplified models of VRFBs, which do not take into account the effects of battery operational parameters on efficiency and lack representation of either capacity fade or restoration mechanisms. Therefore, these models may fail to evaluate accurate optimal dispatching strategies of energy from renewable sources, producing inaccuracies in the prediction of economic revenues from the system. We propose a methodology to include a detailed VRFB modeling in the optimization, building in part on the models in \cite{He2016} and \cite{Jafari2019}, to address the presented problem and to understand how this model performs in terms of optimal results if compared to the standard approach used in literature.
    
    \subsection{Paper contribution}
    The paper addresses the most important gaps found in literature and presents a model containing the following novelty items:
    \begin{itemize}
        \item all the battery losses (cell, stack, auxiliary, inverter) are included in the battery characterization, and the efficiencies are modeled as a function of battery power and state of charge;
        \item efficiency functions and constraints are linearized and convexified, where possible, to obtain a linear problem with a low number of binary variables and a fast solving time;
        \item reversible degradation effects, caused by two different degradation mechanisms, which result in a linear capacity fade with the actual number of cycles of the battery, are included;
        \item two different types of periodic electrolyte maintenance events, needed to maintain the accessible battery capacity above the desired value, are introduced.
    \end{itemize}
    The model is applied to two case studies, and the VRFB stores the energy produced with a renewable energy plant (wind or PV), grid-connected.
    The algorithm decides the optimal operating strategy to maximize the expected profit from the energy exchange with the grid.
    
    The results of this detailed model are then compared to those obtained with simpler models, which do not take into account the capacity fade and consider a constant battery efficiency.
    It is ascertained that having variable efficiencies results in a reduction of the ideal operating range of the battery when compared to the constant efficiencies case, in which the operation at different $SoC$ or at different powers does not penalize the battery operation.
    Furthermore, the modeled capacity fade caused by the previous cycling in the detailed algorithm diminishes the battery revenue, thanks to a decreased accessible energy, and limits the battery usage during the scheduled maintenance period.
    
    The rest of the paper is structured as follows.
    Section \ref{method} describes the methodology and the problem formulation, Section \ref{results} presents the results and discussion, while summary and conclusions are presented in Section \ref{conc}.
    
    \section{Method} \label{method}
    
    This section introduces the problem statement and hypotheses for the model, defining the framework and the case studies to which the problem is applied.
    It also describes the battery and energy system characteristics and defines all the equations and parameters involved in the model.
    
    \subsection{Problem statement} \label{probst}
    The optimization problem, solved separately for each day, during a whole year, can be stated as follows. Given:
    \begin{itemize}[label = -]
        \item size, performance, and state of charge limits of the vanadium redox flow battery system;
        \item renewable energy production curve;
        \item time-dependent demand of electricity;
        \item time-dependent price of selling and purchasing electricity;
    \end{itemize}
    
    determine for each period i of a time horizon T:
    
    \begin{itemize}[label = -]
        \item the battery charging and discharging power;
        \item the state of charge;
        \item the power exchanged with the grid;
    \end{itemize}
    
    to maximize the revenue associated with selling and purchasing electricity from the grid, ensuring that any energy demand is always satisfied.
    
    \subsection{Model formulation} \label{model}
    
    The original optimization problem is a mixed integer and nonlinear. It has a linear objective function, a set of linear constraints, and a set of nonlinear variables and constraints.
    To reduce the computational time and to find a global optimum, the problem is transformed from MINLP (Mixed Integer Non-Linear Program) to MILP with linearization methods.
    
    The hypotheses and modeling strategies are the following:
    
    \begin{itemize}
        \item battery efficiency curves, comprising all system losses, are adapted from empirical data from existing literature \cite{He2016}, inverter losses are also included;
        \item all the nonlinear functions and nonlinear constraints are linearized using different modeling strategies to have a fast linear program with a guaranteed solution of optimality;
        \item the non-convexity of the efficiencies is addressed, and efficiencies are transformed into power functions to reduce the number of required binary variables, therefore reducing the computational burden;
        \item the model requires at least two binary variables (one intrinsic binary, $k_b$ and one additional binary $k_{onoff}$ from the convexification of efficiencies);
        \item due to the low solving time, the problem can be solved with a day-by-day optimization, for a whole year, without the need, for example, for heuristic decomposition techniques such as investigating only a few significant weeks during the annual period
        \item a capacity fade model, from \cite{Rodby2020} is added to the constraints, and the day-by-day optimization allows to solve the additional non-linearity caused by the variable accessible capacity; this works by accounting for historical battery cycling and with daily adjustment of battery capacity;
        \item two types of electrolyte maintenance effects are evaluated both from a technical and economical point of view.
    \end{itemize}
    
    The optimization problem, which is solved with a one-hour resolution time, is formulated based on the hypothesis of having perfect forecasts of energy prices, demand, and renewable energy production for the next 24 hours. Each day is solved separately.
    During each day, the battery charge and discharge phases are optimized, considering the energy prices and the energy production.
    At the end of each day, the daily capacity fade is calculated, and it is passed on the next day.
    
    Within the paper is assumed that the energy efficiency of the battery is influenced both by the charging/discharging power and the state of charge, while the effects of temperature and electrolyte flow are disregarded.
    It is also assumed that the capacity fade of the battery is a linear function of the number of charge/discharge cycles, and the capacity fade rate is constant, thus independent of the current and depth of discharge of the electrolytes \cite{Rodby2020}.
    It is assumed that the efficiency and the degradation effects of the battery unit used to test the optimization model (5 kW/20 kWh Prudent VRFB \cite{Nguyen2014}) are constant for different battery sizes.
    A storage system of any size is made up of many battery units, operated together as an aggregate, maintaining the same performance as a single battery unit.
    
    \subsection{Parameters and variables description}
    
    To describe the optimization model, we introduce the parameters that are known in each period of the day and the variables that are subject to optimization.
    Then we describe step by step the mathematical model by presenting the objective function and all the constraints.
    
    The index \textit{i} indicates a variable or parameter corresponding to the i-th period of the day. The index is chosen from the set T of hours of the day T = \{1,..,24\}.
    The optimization is then carried out for each separate day d for the entire year so that the total period is D = \{1,..,365\}.
    
    \subsubsection{Parameters}
    The problem relies on the following parameters, constant or time-varying:
    \begin{itemize}
        \item $\hat{P}_{b,nom}$  [kW] - Rated battery power;
        \item $\hat{E}_{b,nom}$  [kWh] - Rated battery energy;
        \item $\hat{P}_{res}(i,d)$ [kW] - Power generated by the renewable plant, with i $\in$ T and d $\in$ D;
        \item $\hat{P}_{g,max}$  [kW] - Maximum power that can be exchanged with the grid;
        \item $\hat{P}_{dem}(i,d)$ [kW] - Demanded power by the energy community users, with i $\in$ T and d $\in$ D;
        \item $\hat{P}_{qc}(k)$ \& $\hat{P}_{qd}(k)$ [kW] - Internal charging/discharging battery powers calculated for different linearization sampling points;
        \item $\hat{B}M$ [-] - Big value, enough to deactivate the constraints;
        \item $\hat{SoC}_{min}$ \& $\hat{SoC}_{max}$ [-] - Minimum and maximum state of charge of the battery;
        \item $\hat{SoC}_0$ [-] - Initial state of charge of the battery at the beginning of each day;
        \item $\hat{cap}_{lim}$ [-] - Lower capacity limit of the battery;
        \item $\hat{r}_{ED}$ [\%/cycle] - Electrolyte Decay rate, expressed as \% of rated capacity per complete charge/discharge cycle;
        \item $\hat{R}_{fade}$ [\%/cycle] - Capacity fade rate, expressed as \% of rated capacity per complete charge/discharge cycle;
        \item $\hat{c}_{s}(i,d)$ [€/kWh] -  Price of sold electricity to the grid, with i $\in$ T and d $\in$ D;
        \item $\hat{c}_{p}(i,d)$ [€/kWh] -  Price of purchased electricity from the grid, with i $\in$ T and d $\in$ D.
    \end{itemize}
    
    \subsubsection{Variables}
    The problem makes use of the following decision variables - which are either continuous or binary - leading to a mixed integer program:
    \begin{itemize}
        \item $P_{bc}(i,d)$ $\in \hspace{1mm} \mathbb{R}$ [kW] - 
        Battery charging power, without losses, with i $\in$ T and d $\in$ D;
        \item $P_{bd}(i,d)$ $\in \hspace{1mm} \mathbb{R}$ [kW] - Battery discharging power, with losses, with i $\in$ T and d $\in$ D;
        \item $P_{c,in}(i,d)$ $\in \hspace{1mm} \mathbb{R}$ [kW] - Internal battery charging power, with losses, with i $\in$ T and d $\in$ D;
        \item $P_{d,in}(i,d)$ $\in \hspace{1mm} \mathbb{R}$ [kW] - Internal battery discharging power, without losses, with i $\in$ T and d $\in$ D;
        \item $P_{g,s}(i,d)$ $\in \hspace{1mm} \mathbb{R}$ [kW] - Power sold to the grid, with i $\in$ T and d $\in$ D;
        \item $P_{g,p}(i,d)$ $\in \hspace{1mm} \mathbb{R}$ [kW] - Power purchased from the grid,  with i $\in$ T and d $\in$ D;
        \item $SoC(i,d)$ $\in \hspace{1mm} \mathbb{R}$ [-] - Battery state of charge, with i $\in$ T and d $\in$ D;
        \item $\eta_c(i,d)$ \& $\eta_d(i,d)$ $\in \hspace{1mm} \mathbb{R}$ [-] - Battery charging and discharging efficiencies, with i $\in$ T and d $\in$ D;
        \item $k_b(i,d)$ $\in$ \{0,1\} - Binary variable for charging (1) or discharging (0) the battery, with i $\in$ T and d $\in$ D;
        \item $k_{onoff}(i,d)$ $\in$ \{0,1\} - Binary variable setting the battery state on(1) or off(0), with i $\in$ T and d $\in$ D;
        \item $k_{id}(i,d)$ $\in$ \{0,1\} - Binary variable for purchasing (1) or selling (0) power to the grid, with i $\in$ T and d $\in$ D;
    \end{itemize}
    
    \subsection{Objective function}
    
    The optimization problem consists of the maximization of the following objective function, which is the daily net revenue coming from the energy sold to the grid minus the expenses for energy purchasing:
    
    \begin{equation}\label{FO}
        Rev(d) = \tau \cdot \sum_{i=1}^{T} \hat{c}_{s}(i,d) \cdot P_{g,s}(i,d) - \tau \cdot \sum_{i=1}^{T} \hat{c}_{p}(i,d) \cdot P_{g,p}(i,d)
    \end{equation}
    
    Where $\tau$ is the period over which the powers are assumed to be constant. In this work, it is always assumed that $\tau = 1$ h, therefore T = 24.
    $\hat{c}_s$ and $\hat{c}_p$ are the cost of electricity sold and purchased to and from the grid.
    $P_{g,s}$ and $P_{g,p}$ are the power sold and purchased to and from the grid.
    
    \subsection{Constraints} \label{constraints}
    
    \subsubsection{Power constraints}
    
    The relationship between the battery, grid, load, and renewable source is expressed through the curtailment power, which is the sum of all the power terms:
    
    \begin{equation} \label{curtailment}
        P_{curt}(i,d) = \hat{P}_{res}(i,d) - P_b(i,d) - P_{g,s}(i,d) + P_{g,p}(i,d) - \hat{P}_{dem}(i,d) \geq 0
    \end{equation}
    
    The power constraint (\ref{curtailment}) is relaxed with an inequality to increase solving speed. The algorithm tends to minimize $P_{curt}$ for the optimization.
    
    The battery power $P_{b}$ is expressed as follows:
    
    \begin{equation}
        P_b(i,d) = P_{bc}(i,d) - P_{bd}(i,d)
    \end{equation}
    
    Grid power is subjected to the following constraints to avoid simultaneously selling and purchasing electricity from the grid:
    
    \begin{equation}\label{pg1}
        0 \leq P_{g,p}(i,d) \leq k_{id}(i,d) \cdot \hat{P}_{g,max}
    \end{equation}
    \begin{equation}\label{pg2}
        0 \leq P_{g,s}(i,d) \leq \left(1 - k_{id}(i,d)\right)\cdot \hat{P}_{g,max}
    \end{equation}
    
    To guarantee that the battery can not simultaneously charge and discharge, the following constraints are introduced:
    
    \begin{equation}\label{v1}
        0 \leq P_{bc}(i,d) \leq k_b(i,d) \cdot \hat{P}_{b,nom}
    \end{equation}
    \begin{equation}\label{v2}
        0 \leq P_{bd}(i,d) \leq (1 - k_b(i,d)) \cdot \hat{P}_{b,nom}
    \end{equation}
    \begin{equation}\label{v3}
        -k_b(i,d) \cdot \text{max}(\hat{P}_{qc}(k)) \leq P_{c,in}(i,d) \leq k_b(i,d) \cdot \text{max}(\hat{P}_{qc}(k))
    \end{equation}
    \begin{equation}\label{v4}
        0 \leq P_{d,in}(i,d) \leq (1 - k_b(i,d)) \cdot \text{max}(\hat{P}_{qd}(k))
    \end{equation}
    
    Additionally, there are constraints for the on/off state of the battery:
    
    \begin{equation}\label{v5}
        P_{bc}(i,d) \leq k_{onoff}(i,d) \cdot \hat{P}_{b,nom}
    \end{equation}
    \begin{equation}\label{v6}
        P_{bd}(i,d) \leq k_{onoff}(i,d) \cdot \hat{P}_{b,nom}
    \end{equation}
    \begin{equation}\label{v7}
        P_{c,in}(i,d) \leq k_{onoff}(i,d) \cdot \text{max}(\hat{P}_{qc}(k))
    \end{equation}
    \begin{equation}\label{v8}
        P_{d,in}(i,d) \leq k_{onoff}(i,d) \cdot \text{max}(\hat{P}_{qd}(k))
    \end{equation}
    
    \subsubsection{State of charge and state of energy constraints}
    
    The battery state of charge ($SoC$) is defined as the mean energy inside the battery during the time interval \textit{i} (SoE, kWh) divided by the maximum energy that can be stored in the battery ($E_m(d-1)$, kWh), following a similar approach to the one in \cite{He2016}:
    
    \begin{equation}
        SoC(i,d) = \frac{SoE(i,d) + SoE(i-1,d)}{2\cdot E_m(d-1)}
    \end{equation}
    
    During the day, the battery experiences side reactions in the cell, which result in a capacity fade, and a new value of capacity is calculated at the end of each day.
    This is done to solve the degradation nonlinearity, optimizing each day separately and then passing the results on to the next one. $E_m(d-1)$ is the accessible energy that the battery can store at the end of day \textit{d-1}, or equivalently at the beginning of day \textit{d}.
    
    Note that some authors \cite{Zugschwert2021, He2016} define the battery state of charge ($SoC$) in terms of the maximum accessible capacity, as done here, while others \cite{Jafari2019, Gong2019}, define $SoC$ always in terms of the rated battery capacity.
    In this work, we decided to define $SoC$ in terms of the accessible capacity because this definition is considered to be more adherent to the real status of the battery at the expense of introducing a nonlinear term to the model.
    In this way, in fact, $SoC$ is defined as the ratio of two functions of time and is, thus, a nonlinear term.
    The nonlinearity is solved by optimizing each day (d) separately from the others and defining $SoC$ in terms of the maximum accessible energy from the previous day (d-1).
    
    $SoC$ and $SoE$ are bound between the maximum and minimum levels:
    
    \begin{equation}\label{vsoc}
        \hat{SoC}_{min} \leq SoC(i,d) \leq \hat{SoC}_{max}
    \end{equation}
    \begin{equation}\label{vsoe}
        0 \leq SoE(i,d) \leq E_m(d-1)
    \end{equation}
    
    $\hat{SoC}_{min}$ and $\hat{SoC}_{max}$ are set to be 10\% and 90\%, respectively. These boundaries define a commonly used $SoC$ range for VRFBs, set with the purpose of avoiding extreme operating conditions, where charging could occur at lower efficiency, or where the magnitude of hydrogen evolution becomes significant (i.e. above 90\%) \cite{Tang2011}.
    
    Also, SoE has to be bounded with the following constraint every day:
    
    \begin{equation}\label{vsoend}
        SoE(i,d)\Bigr|_{\substack{i = 24}} = \hat{SoE}_0
    \end{equation}
    
    Where the initial energy contained in the battery is:
    
    \begin{equation}
        \hat{SoE}_0 = \hat{SoC}_0 \cdot \hat{E}_{b,nom}
    \end{equation}
    
    \subsubsection{Battery efficiency linearization}
    
\begin{figure}[ht]
        \centering
        \subfigure[]{%
        \includegraphics[scale = 0.6]{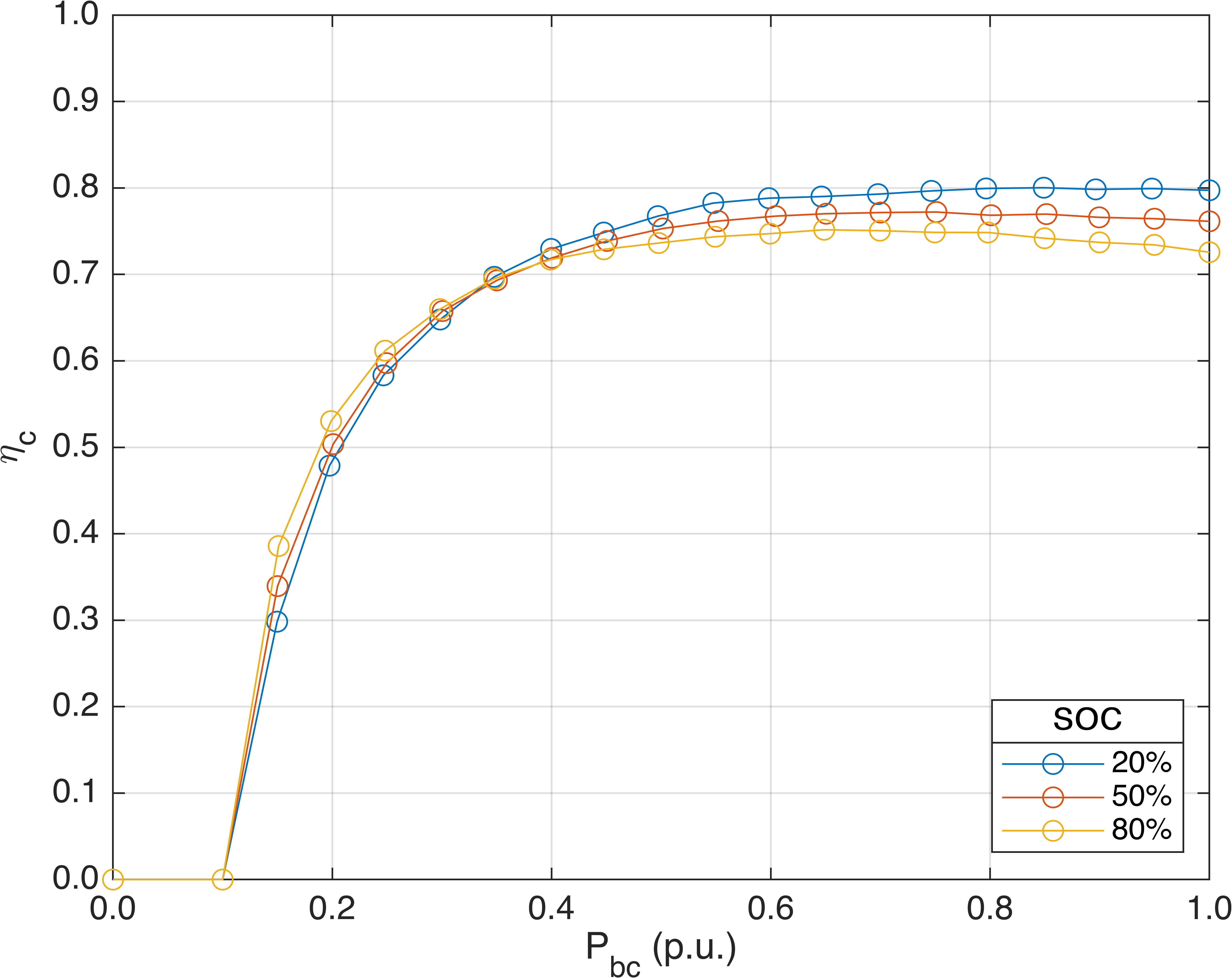}
        }
        \quad
        \subfigure[]{%
        \includegraphics[scale = 0.6]{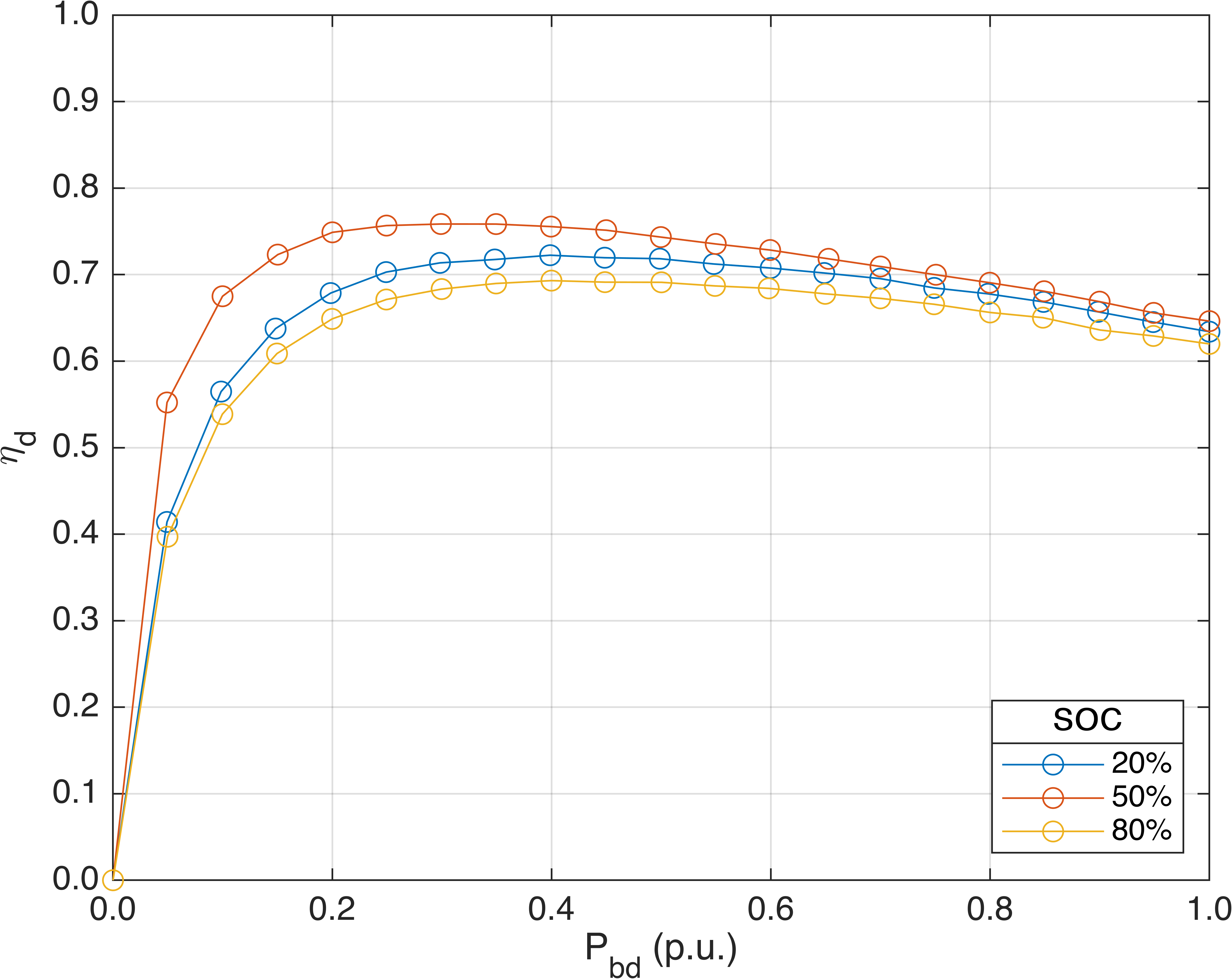}
        }
        \caption{Charging (a) and discharging (b) efficiencies for three different $SoC$ levels and dimensionless power values. Data source: \cite{He2016}}
        \label{efficienze}
    \end{figure}

\begin{figure}[hbt]
        \centering
        \includegraphics[scale = 0.5]{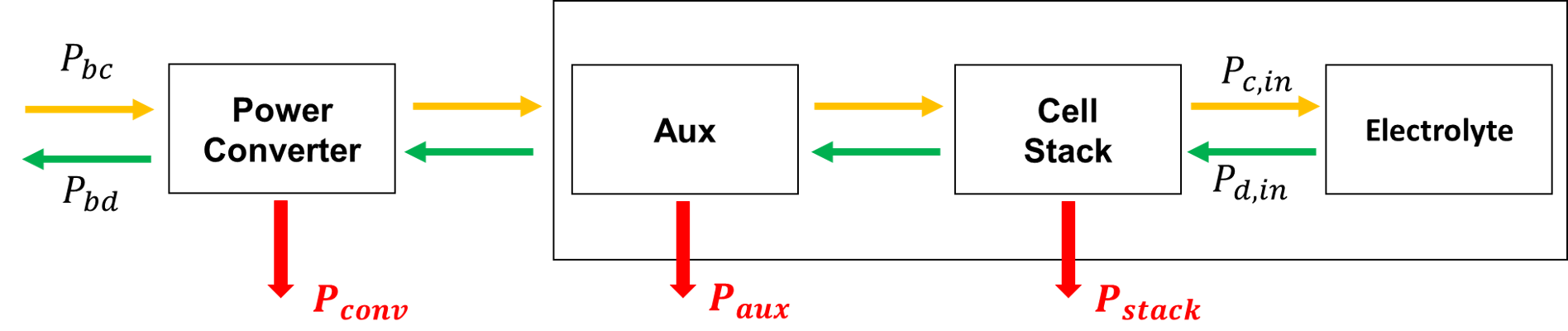}
        \caption{Battery simplified scheme with powers and power losses}
        \label{schema_tot}
    \end{figure}

    The battery efficiency is modeled as a function of both charging/discharging power and $SoC$.
    The data shown in Fig. \ref{efficienze} represent the charge and discharge global efficiency curves for the 5kW/20kWh Prudent VRFB at different $SoC$ levels (20\%, 50\%, and 80\%) and different charge and discharge power values, expressed in per-units (p.u.), relatively to the rated power.
    These energy efficiencies include all the internal cell stack losses, plus the auxiliary systems losses (pumps with their dedicated electric motor and inverter, and battery controls), including the inverter losses (see Figure \ref{schema_tot} and relative description for a simplified scheme of losses location).
    The experimental data have been conducted at a constant temperature of 25°C, at the maximum flow rate allowed by the pumps.
    The HVAC (Heat, Ventilation, and Air Conditioning) system is usually deployed to maintain a constant temperature in the battery. The energy losses in the HVAC system are not considered in this model because they have little impact on the optimal solution \cite{He2016} and will depend on the temperature of the installation site.
    
    As seen in Fig. \ref{efficienze}, the efficiencies are nonlinear with $SoC$ and power and for both cases, at low power, the efficiency rapidly decreases, mainly because the auxiliary pumps absorb much of the provided energy, causing higher energy losses.
    
    The state of energy of the battery depends on the energy efficiencies as follows:
    
    \begin{equation} \label{SoE_eq}
        SoE(i,d) = \hat{SoE}_0 + \tau \cdot \sum_{n=1}^{i} \eta_{c}(SoC,P_{bc})  \cdot P_{bc}(n,d) - \tau  \cdot \sum_{n=1}^{i} \frac{1}{\eta_{d}(SoC,P_{bd})} \cdot P_{bd}(n,d)
    \end{equation}
    
   The latter equation is a source of nonlinearity in the model because it encompasses terms of multiplication and division between decision variables of the problem ($P_{bc}$ and $\eta_c$ or $P_{bd}$ and $\eta_d$).
   To solve this nonlinearity, $SoE$ is then defined as follows:
   
    \begin{equation} \label{SoE_eq2}
        SoE(i,d) = \hat{SoE}_0 + \tau \cdot \sum_{n=1}^{i} P_{c,in}(i,d) - \tau \cdot \sum_{n=1}^{i} P_{d,in}(i,d)
    \end{equation}
    
    Where:
    
    \begin{equation} \label{charging_power}
        P_{c,in} = \eta_c \cdot P_{bc} = P_{bc} - P_{c,losses}
    \end{equation}
    \begin{equation} \label{discharging_power}
        P_{d,in} = \frac{1}{\eta_d} \cdot P_{bd} = P_{bd} + P_{d,losses}
    \end{equation}
    
    Equations (\ref{charging_power}) and (\ref{discharging_power}) delineate the relationship between the internal powers and external powers of the battery, as defined in this work, through battery power losses.
    $P_{bc}$ is the external charging power, as seen from the external battery connection, $P_{c,in}$ is the internal charging power, as calculated from ideal electrolyte properties, without considering losses $P_{c,losses}$, and the relative effect applies to the discharging powers.
    The power losses terms ($P_{c,losses}$, $P_{d,losses}$) and the efficiencies ($\eta_c$, $\eta_d$) encompass all the losses in the battery.
    These losses are located, as shown in Figure \ref{schema_tot}, not only in the cell stack due to both intrinsic cell losses and parasitic shunt currents but also in the auxiliary system (which comprehends circulation pumps and other power-control devices) and, lastly, in the inverter or static power converter.
    
\begin{figure}[ht]
        \centering
        \subfigure[]{%
        \includegraphics[scale = 0.6]{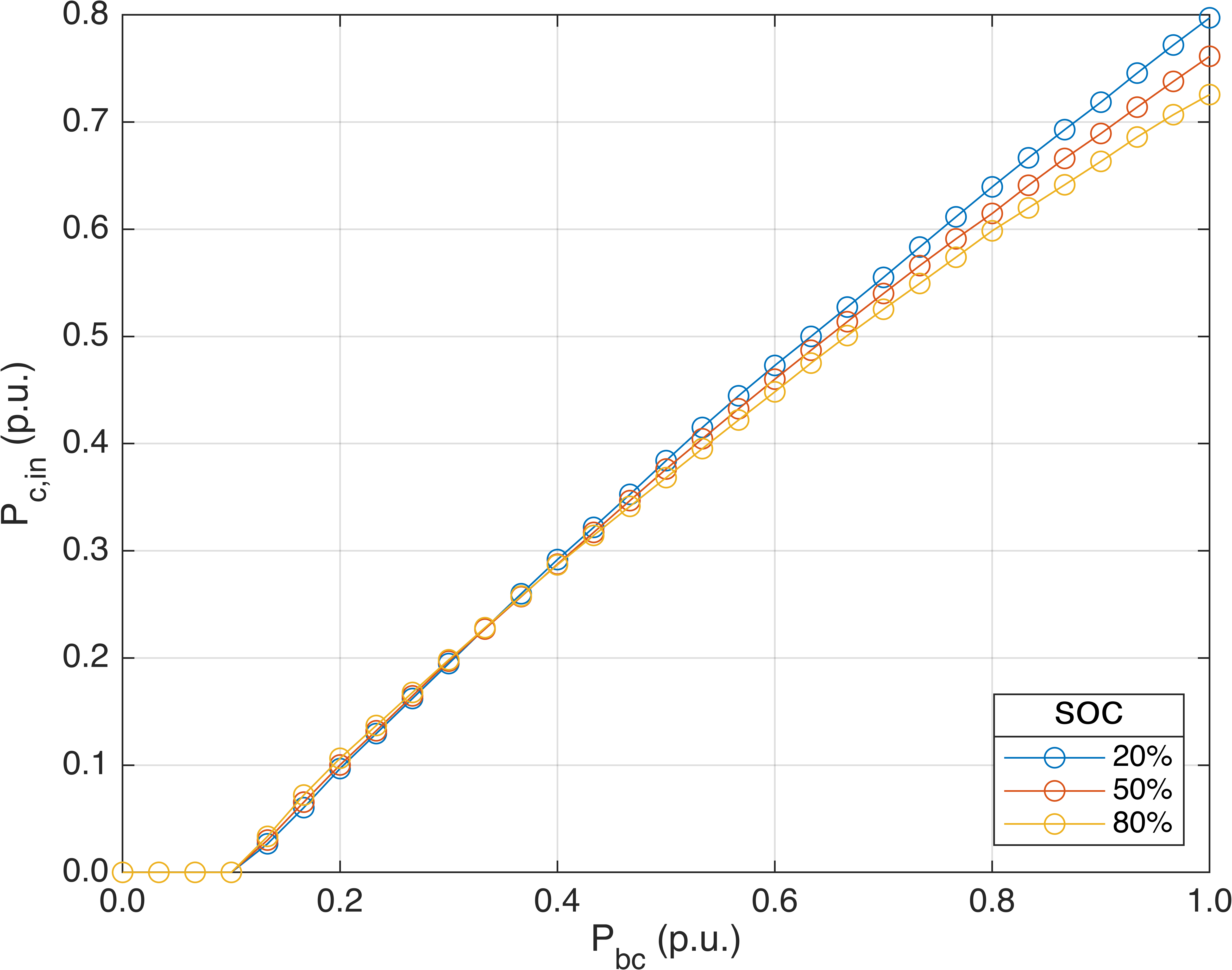}
        }
        \quad
        \subfigure[]{%
        \includegraphics[scale = 0.6]{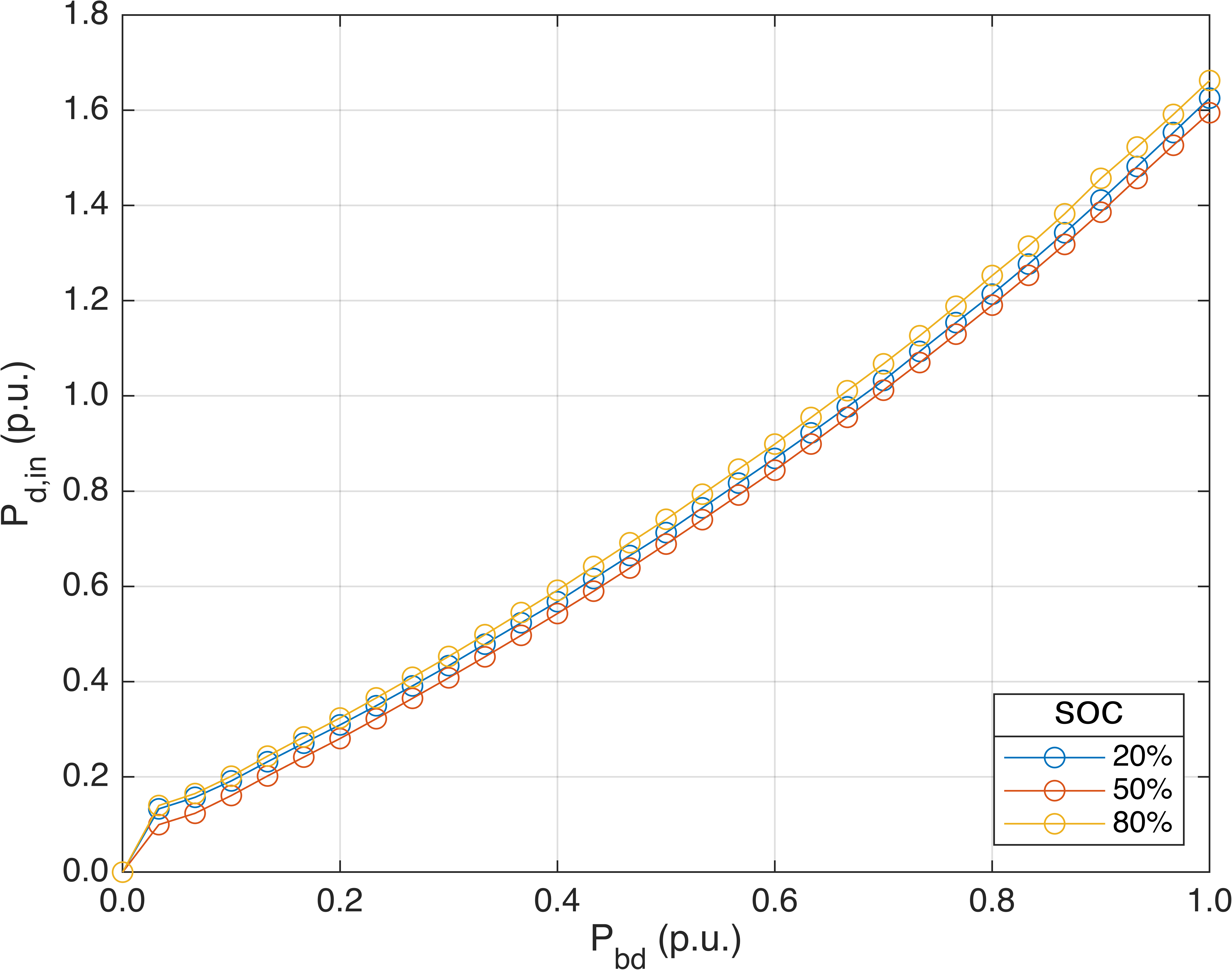}
        }
        \caption{Charge (a) and discharge (b) internal powers for three different $SoC$ levels and dimensionless power values}
        \label{p_internal}
    \end{figure}

    As shown in Figure \ref{p_internal}, $P_{c,in}$ and $P_{d,in}$ are also nonlinear, non-convex functions of dimensionless charging/discharging power (expressed in per-units, relative to the rated battery power) and $SoC$.
    To linearize them, we enclose these internal power functions in convex hulls with the use of planes tangent to the curves to obtain linear variables and constraints while introducing a binary variable as a result of the convexification.
    
\begin{figure}[hbt]
        \centering
        \subfigure[]{%
        \includegraphics[scale = 0.58]{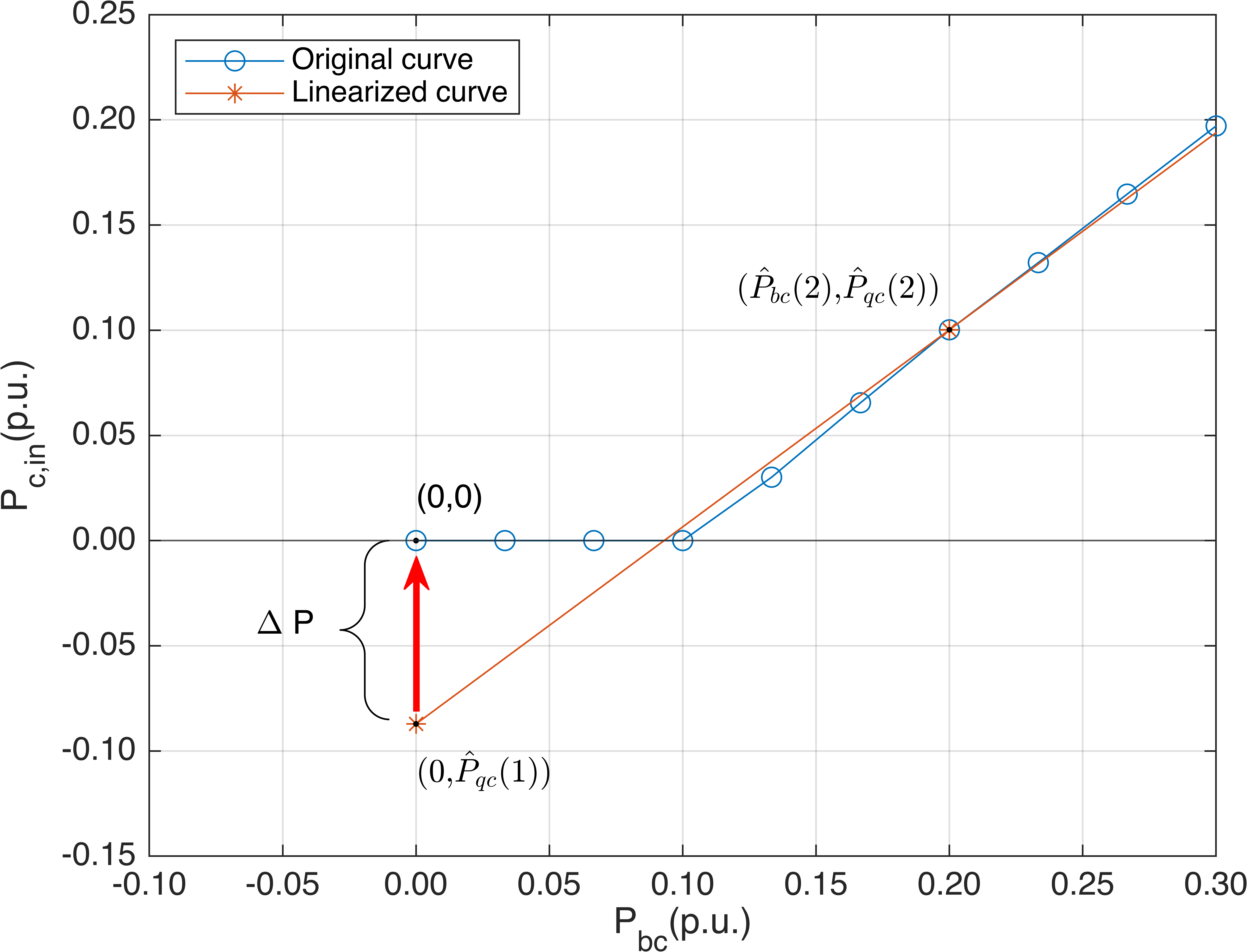}
        }
        \quad
        \subfigure[]{%
        \includegraphics[scale = 0.58]{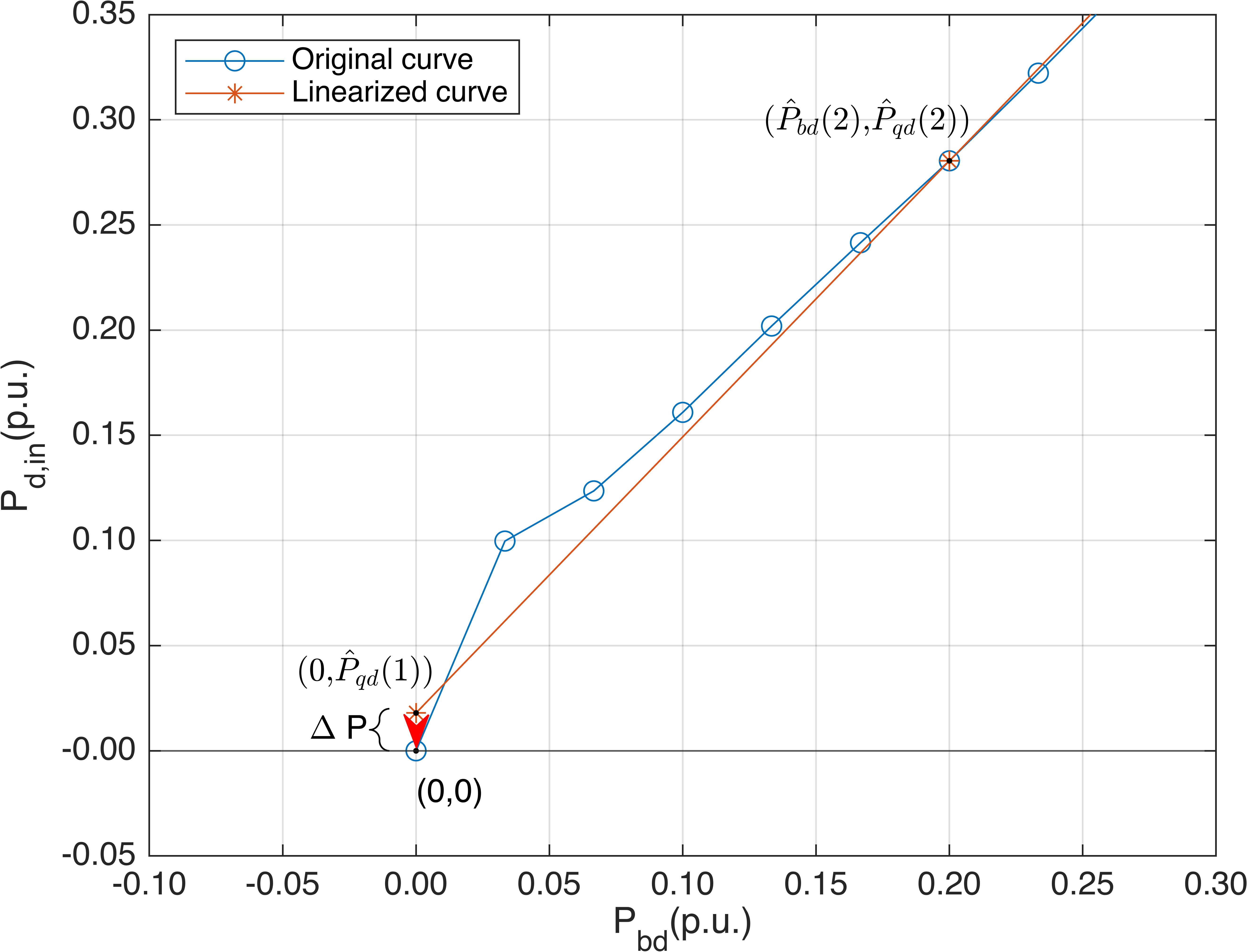}
        }
        \caption{Close-up graph of original and linearized convex functions of $P_{c,in}$ (a) and $P_{d,in}$ (b), for $SoC$ = 50\%}
        \label{linearize}
    \end{figure}

    The two-variable linearization is created by dividing the power domain into $n_{int}$ intervals of equal length, while $SoC$ is interpolated between the three originally known values of 20, 50, and 80\%.
    For each of the three $SoC$ curves, there are $n_{int}+1$ linearization sampling points in regard to the charging and discharging powers. The sampling points are ($\hat{P}_{bc}(k)$, $\hat{P}_{qc}(k)$) and ($\hat{P}_{bd}(k)$, $\hat{P}_{qd}(k)$), where k = 1,...,$3\cdot (n_{int}+1) $. Figure \ref{linearize} shows a close-up of the linearized curve for $SoC = 50\%$, where the first two sampling points are visible.
    
\begin{figure}[ht]
        \centering
        \subfigure[]{%
        \includegraphics[scale = 0.59]{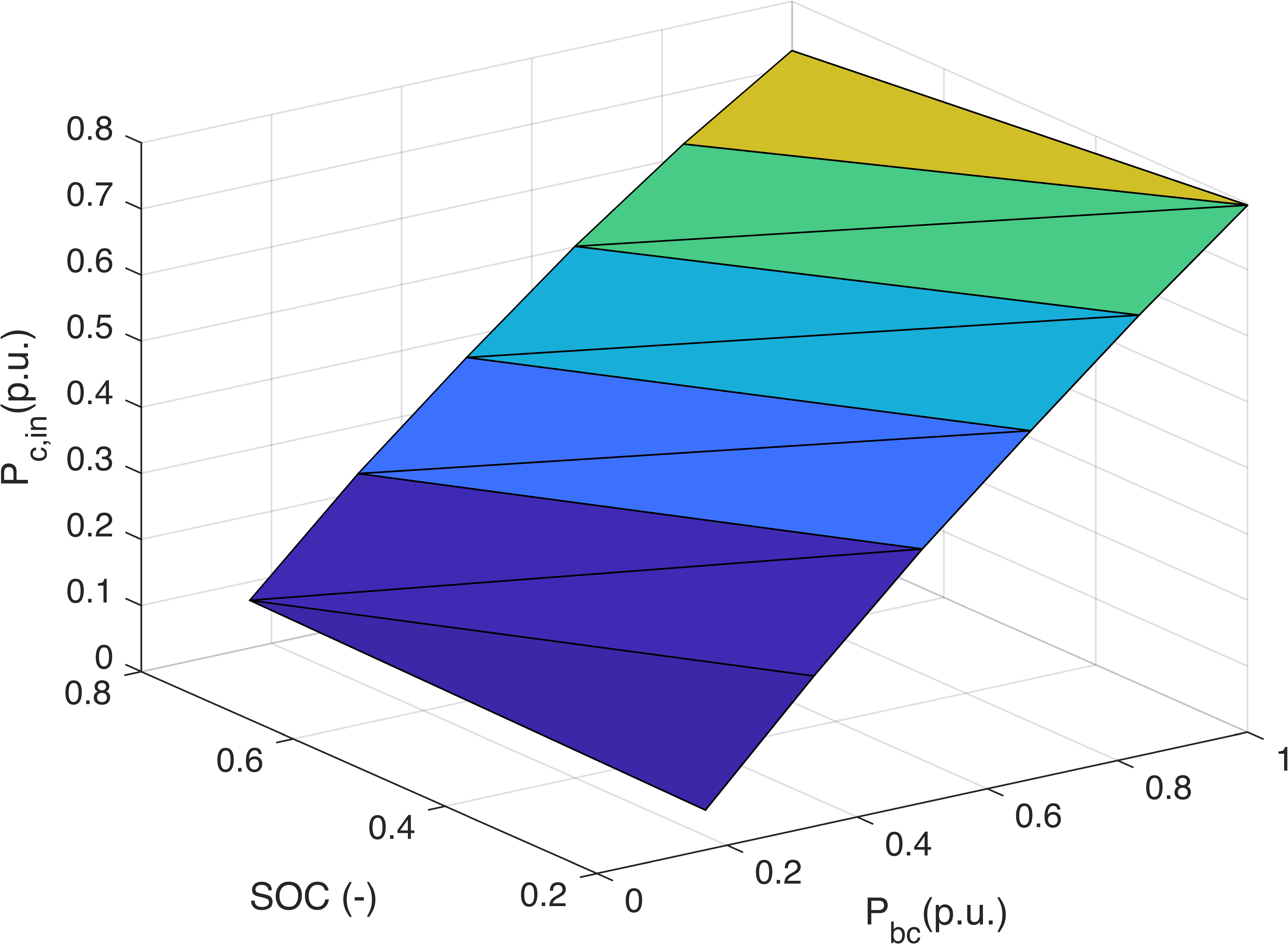}
        }
        \quad
        \subfigure[]{%
        \includegraphics[scale = 0.59]{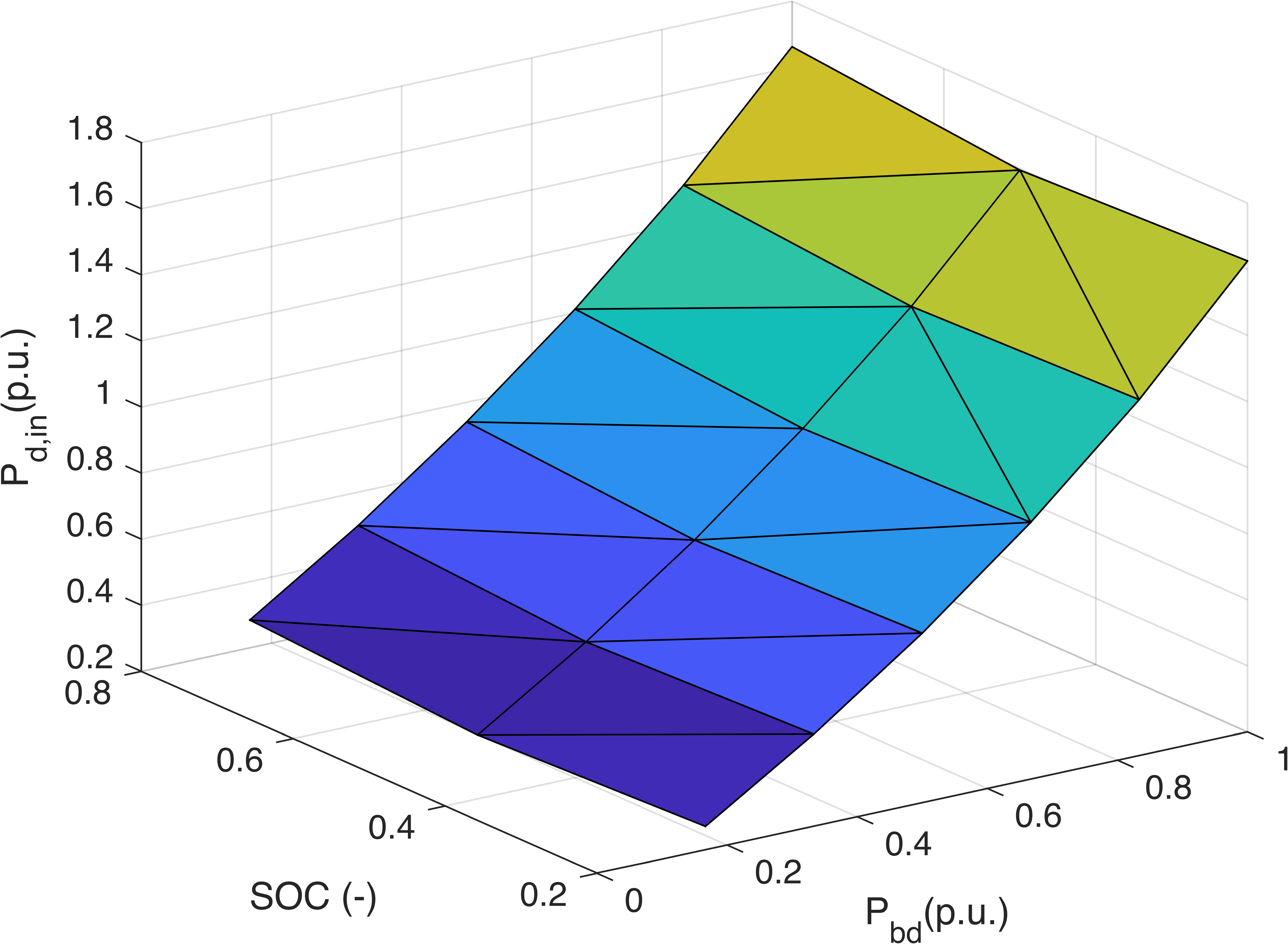}
        }
        \caption{Convex hull selected tangent planes from the linearization of the two-variables functions of internal powers ($P_{c,in}$ (a) and $P_{d,in}$ (b)), with $n_{int}$ = 5}
        \label{convhull}
    \end{figure}

    The tangent planes for the linearization of the internal power curves in Figure \ref{p_internal} are selected among those creating a \textit{convex hull} containing the functions, and can be seen in Figure \ref{convhull}.
    
    The $J_c$ tangent planes to the two-variable charging power surfaces, creating the convex hull of the charging power, are defined by the following equations:
    
    \begin{equation}
        z = \hat{\gamma}_{1,jc} \cdot x+\hat{\gamma}_{2,jc} \cdot y+\hat{\gamma}_{3,jc} \hspace{0.5 cm} \text{with } j_c = 1, ...., J_c
    \end{equation}
    
    The $J_d$ tangent planes to the two-variable discharging power surfaces, creating the convex hull of the discharging power, are defined by the following equations:
    
    \begin{equation}
        z = \hat{\delta}_{1,jd} \cdot x+\hat{\delta}_{2,jd} \cdot y+\hat{\delta}_{3,jd} \hspace{0.5 cm} \text{with } j_d = 1, ...., J_d
    \end{equation}
    
    For every $j_c$ tangent plane to the $P_{c,in}$ curves and every $j_d$ tangent plane to the $P_{d,in}$ curves, the following constraints, which state the lower and upper bounds of the internal power functions, are defined as follows:
    
    \begin{equation} \label{pccons}
        P_{c,in}(i,d) \leq \hat{\gamma}_{1,jc}\cdot P_{bc}(i,d) + \hat{\gamma}_{2,jc}\cdot SoC(i,d) + \hat{\gamma}_{3,jc} + \hat{B}M \cdot (1-k_b(i,d)) + \hat{B}M \cdot (1-k_{onoff}(i,d)) 
    \end{equation}
    \begin{equation} \label{pdcons}
        P_{d,in}(i,d) \geq \hat{\delta}_{1,jd}\cdot P_{bd}(i,d) + \hat{\delta}_{2,jd}\cdot SoC(i,d) + \hat{\delta}_{3,jd} - \hat{B}M \cdot k_b(i,d) - \hat{B}M \cdot (1-k_{onoff}(i,d)) 
    \end{equation}
    
    The optimal solution is achieved when the energy contained inside the battery has the maximum value, which means minimum losses, for a given charging/discharging power value.
    This means that, as seen in Eq. (\ref{SoE_eq2}), the solver will find the optimal solution by trying to maximize $P_{c,in}$ and minimize $P_{d,in}$.
    Therefore, we use Eq. (\ref{pccons}) and (\ref{pdcons}) to set upper and lower bounds to the internal battery powers.
    These constraints are linear: the first set of equations (\ref{pccons}) defines the charging power upper bounds and the second one (\ref{pdcons}) defines the discharging power lower bounds.
    
    In these equations, $\hat{B}M$ is a parameter with a big positive value, used to deactivate the constraints when they are not needed, according to the values of the binary variables $k_b$ and $k_{onoff}$.
    To avoid simultaneous battery charging and discharging, the power constraints use the binary variable $k_b$ (set to 1 during charging and to 0 during discharging), to activate one constraint at a time.
    Furthermore, an on/off binary variable ($k_{onoff}$) needs to be introduced in the model to deactivate the constraints when the battery is turned off and both $P_{c,in}$ and $P_{d,in}$ needs to reach zero.
    While the charging/discharging binary variable $k_b$ is an intrinsic variable of the problem, the on/off variable $k_{onoff}$ needs to be introduced due to the convexification of the battery internal powers, obtained by transforming the non-convex battery efficiencies.
    
    As seen in Figure \ref{p_internal}, $P_{c,in}$ is convex upward (or concave) with respect to $P_{bc}$, only for $P_{bc} > 0.1$, and $P_{d,in}$ is convex downward (or convex) with respect to $P_{bd}$, only for $P_{bd} > 0.03$. $P_{bd}$'s lower limit value comes from an approximation of the battery power curves. As shown in Fig. \ref{efficienze}, when $P_{bd}$ approaches zero, the discharging efficiency decreases rapidly. From Eq. (\ref{discharging_power}) it then follows that when the discharging efficiency is very low (approximately when $P_{bd}$ is below 0.03), $P_{d,in}$ starts to increase and tends to infinity. To avoid indefinite values and physical infeasibility, $P_{d,in}$ value is set to zero when $P_{bd}$ equals zero.
    If $P_{c,in}$ were a convex upward function and $P_{d,in}$ were a convex downward function through all their domain of existence, it would be possible to perform a linear approximation using tangent planes to the two-variable function while guaranteeing that the solver would push the solution as closely as possible to the approximating pieces of curves. This problem would therefore be convex and would not require an additional binary variable ($k_{onoff}$) to be solved.
    Because, instead, $P_{c,in}$ is convex upward only above 0.1 and $P_{d,in}$ is convex downward only above 0.03, we first use convex hulls to convexify the functions, and then introduce binary variables to cope with the error from the convexification.
    The error of the convexification ($\Delta P$), can be seen in Figure \ref{linearize}, which shows the original curve and the linearized curve of the battery powers for a fixed $SoC$ value.
    
    When the battery is turned off ($P_{bc} = P_{bd} = 0$), the linearized curves will assume the values $\hat{P}_{qc}$ and $\hat{P}_{qd}$, which differ from the values of the original curves (zero) of an error $\Delta P$, which means that the battery keeps consuming energy.
    To compensate $\Delta P$ and to set $P_{c,in} = P_{d,in} = 0$, when $P_{bc} = P_{bd} = 0$, we use, for both charging and discharging powers, a $k_{onoff}$ binary variable to activate a constant parameter $\hat{B}M$, added up or subtracted to the linearized functions of the internal powers, as seen in Eq. (\ref{pccons}) and (\ref{pdcons}).
    The role of $\hat{B}M$ is to relax the constraints enough to deactivate them when they are not required by the problem. This feature is enabled by the binary variables, which act as an "if condition".
    In particular $\hat{B}M$ deactivates the following constraints:
    \begin{itemize}
        \item both constraints, when the battery is in off ($k_{onoff} = 0$) and it is not charging, nor discharging;
        \item the charging constraint, when the battery is discharging ($k_b = 0$, $k_{onoff} = 1$);
        \item the discharging constraint, when the battery is charging ($k_b = 1$, $k_{onoff} = 1$).
    \end{itemize}
    
    \text{green}{Note that, when the battery is off (i.e. $k_{onoff}=0$), $k_b$ is free to take any arbitrary value because all the relevant variables are bounded, and the constraints are deactivated by $k_{onoff}$ so that the charging/discharging binary variable $k_b$ does not affect the problem.}
    
    For the algorithm to work, the value of $\hat{B}M$, has to be the following:
    
    \begin{equation}
        \hat{B}M = \Delta P + O(\Delta P);
    \end{equation}
    
    \subsubsection{Capacity fade and recovery} \label{capacity_fade_cons}
    
    As previously stated in the Introduction, the most important degradation effects occurring in the vanadium battery are crossover and undesired oxidation reactions, which cause a reduction in the accessible capacity of the battery.
    Crossover, caused by concentration and pressure differences in the two half-cells, is one of the biggest drivers of the capacity fade.
    In the VRFB, the crossover does not cause permanent cross-contamination between the two sides of the cell because of the symmetric chemistry of the battery. Thanks to this property, the crossover effects can be reversed using the rebalancing method, which consists of periodic electrolyte mixing followed by complete recharging \cite{Rodby2020, Poli2020}.
    The undesired oxidation reactions are less detrimental to capacity than the crossover effect, but the process to reverse their effects (servicing) is more expensive and more difficult to perform than the process to reverse crossover effects (rebalancing). Undesired side reactions lead to the oxidation of vanadium ions from contact with air, vanadium gassing or loss of electrolyte stability due to impurities, and hydrogen evolution on the negative electrode, where the electrons react with hydrogen ions instead of reacting with vanadium ions \cite{Yuan2019}.
    These secondary reactions are responsible for an oxidative imbalance between the two sides of the cell, resulting in a difference between the quantities of reduced and oxidized species in the electrolytes, which consequently induces a capacity fade.
    There are also other types of irreversible degradation of the battery components, which cause efficiency and performance losses over time, and have both terms of the calendar and cyclic degradation \cite{Yuan2019}, but they are very limited and have little influence on the capacity fade; \cite{Rodby2020}, therefore, they will be neglected here.
    
    In this model, the capacity fade process and recovery for the vanadium battery can be explained through the following hypothesis, as assumed by Rodby et al. \cite{Rodby2020}:
    \begin{itemize}
        \item The rated energy of the battery, defined as the product between rated power and discharge duration, is 100\% at the beginning of its operating life.
        \item The capacity fade depends linearly on the number of charge/discharge cycles, while calendar degradation of other components on the performance is neglected. The capacity fade rate is constant throughout the whole battery life.
        \item The effects of the current density, the depth of discharge, and the average oxidation state on the accessible capacity of a battery are neglected.
        \item After fading to a lower capacity limit ($\hat{cap}_{lim}$), the capacity is restored by a complete rebalancing event to a new value of capacity: $cap_{max}$.
        \item The new capacity equals the maximum battery capacity at the last rebalancing event, minus the fade from oxidative imbalance degradation.
        \item Capacity fade and rebalancing events occur until the $cap_{max}$ itself decreases to the $\hat{cap}_{lim}$ when a servicing event is performed. After a servicing event, all the faded capacity is restored to the initial value.
        \item During the battery life, these steps are repeated continuously.
    \end{itemize}
    Note that battery capacity and energy are proportional to each other, as in this optimization work, the potential of the battery is considered constant; therefore, talking about a relative capacity fade means talking about an equal amount of relative energy fade.
    
    As seen in section \ref{VRFB}, capacity fade is caused by crossover (concentration imbalance) and secondary oxidative reactions (oxidative imbalance):
    \begin{itemize}
        \item $\hat{r}_{CR}$ - crossover fade rate;
        \item $\hat{r}_{ED}$ - electrolyte decay rate due to oxidative imbalance;
        \item $\hat{R}_{fade}$ - total capacity fade rate;
    \end{itemize}
    \begin{equation}
        \hat{R}_{fade} = \hat{r}_{CR} + \hat{r}_{ED}
    \end{equation}
    According to a VRFB review from Rodby et al. \cite{Rodby2020}, the total $\hat{R}_{fade}$ stands between 0.067\% and 1.3\%, for different types of batteries. This value is strongly influenced by technology.
    For the current work, the total fade rate is set at $\hat{R}_{fade} = 0.442 \%$, while the electrolyte decay rate is evaluated at $\hat{r}_{ED} = 0.055\%$. The electrolyte decay rate depends upon the electrolyte type and the operating conditions.
    
    \begin{figure}[hbt]
        \centering
        \includegraphics[width = 0.8\textwidth]{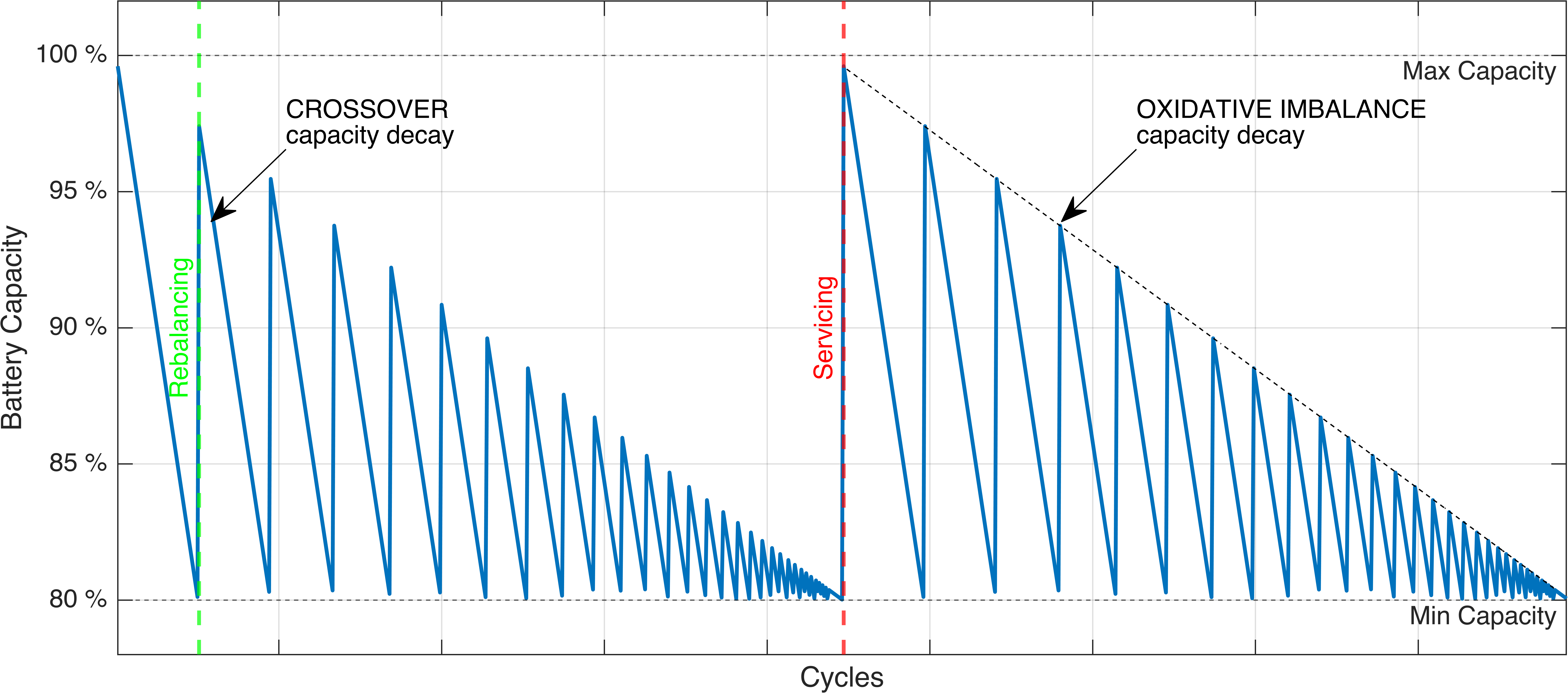}
        \caption{Qualitative simulation of relative battery capacity to the rated value, as a function of the number of cycles and the effects of rebalancing and servicing \cite{Rodby2020}}
        \label{cap_decay}
    \end{figure}
    
    As seen in the example in Figure \ref{cap_decay}, the lower capacity limit ($\hat{cap}_{lim}$) is set at 80\%, because as demonstrated by Rodby et al. \cite{Rodby2020}, the optimal value of the capacity limit resides around 80-85\% for VRFBs with the considered characteristics.
    Maintenance events are represented by vertical lines, which partially or totally restore the battery capacity to maintain it above the lower limit.
    Given that the crossover effect, mainly responsible for the concentration imbalance, is of a larger magnitude compared to the oxidation effects, the rebalancing (first event signaled by a green line) occurs more frequently than the servicing (signaled by a red line).
    
    The equation that determines the accessible energy of the battery for each day d is the following:
    
    \begin{equation}
        E_m(d) = \hat{E}_{b,nom} \cdot f_{cap}(d)
    \end{equation}
    
    Where $f_{cap}(d)$ is the fraction of the accessible capacity ad the beginning of day d.
    $f_{cap}(d)$ is calculated using the number of cycles passed since the last rebalancing ($n_{cyc}^R(d)$) or servicing ($n_{cyc}^R(d)$) event:
    
    \begin{equation}
        f_{cap}(d) = cap_{max}(d) - \hat{R}_{fade}\cdot n_{cyc}^R(d)
    \end{equation}
    \begin{equation}
        cap_{max}(d) = 1 - \hat{r}_{ED}\cdot n_{cyc}^S(d)
    \end{equation}
    
    To evaluate the number of cycles passed since the last maintenance event ($n_{cyc}^R$ e $n_{cyc}^S$), an algorithm counts the energy stored inside the battery:
    
    \begin{equation}
        n_{cyc}^{R/S}(d) = \begin{cases}
      0 & \text{immediately after \textit{rebalancing} (R) or \textit{servicing} (S)}\\
      n_{cyc}^{R/S}(d-1) + N_{cyc,d}(d) & \text{in between R or S events}\\
        \end{cases} 
    \end{equation}
    
    \begin{equation} \label{cycling}
        N_{cyc,d}(d) = \frac{\tau \cdot \sum_{i=1}^{24} P_{c,in}(i,d)}{\hat{E}_{b,nom}}
    \end{equation}
    
    Where $N_{cyc,d}(d)$ is the number of partial charge/discharge cycles during the day d or the amount of energy charged into the battery for that day divided by the rated capacity of the battery. On a daily basis, the energy charged into the battery equals the energy discharged from the battery, because the self-discharge is neglected, as it is generally considered very low for vanadium flow batteries \cite{Diaz2012}.
    Note that partial charge and discharge cycles partially degrade the battery capacity.
    
    Maintenance events are executed when the accessible capacity reaches the lower capacity limit $\hat{cap}_{lim}$, which is the lowest fraction of allowed accessible capacity for the battery:
    \begin{itemize}[label = -]
        \item if $f_{cap}(d) \leq \hat{cap}_{lim}$ a \textit{rebalancing} is needed;
        \item if $cap_{max}(d) \leq \hat{cap}_{lim}$ a \textit{servicing} is needed.
    \end{itemize}
    
    To model the capacity fade and, consequently, the maintenance events, additional rebalancing constraints are needed:
    
    \begin{equation} \label{reb_cons1}
        SoC(i,d)\Bigr|_{\substack{i = \hat{i}_{reb}, d = d_{reb}}} = \hat{SoC}_{max}
    \end{equation}
    \begin{equation} \label{reb_cons2}
        P_{bd}(i,d) \Bigr|_{\substack{i = 1:\hat{i}_{reb}, d = d_{reb}}} = 0
    \end{equation}
    
    Eq. (\ref{reb_cons1}) and (\ref{reb_cons2}) are valid only during the rebalancing day $d_{reb}$ and during the rebalancing event, which starts at the beginning of the day ($i=1$) and ends at $\hat{i}_{reb}$.
    These constraints define the state of the battery during and after every rebalancing event.
    The rebalancing takes place in two steps: the first is the electrolyte mixing (which is assumed to be fast and whose occurrence time is not modeled), and the second is the electrolyte recharging from a state of complete homogeneity between the two tanks and the two electrolytes, up to the maximum possible state of charge, as stated by Eq. (\ref{reb_cons1}).
    During the recharging, which is assumed to last from the first hour of the rebalancing day, up to $\hat{i}_{reb}$, the battery cannot be used, and it is allowed to recharge with the necessary energy for complete recharging of the electrolyte. During this time, the discharging power of the battery is, therefore, set to zero by Eq. (\ref{reb_cons2}).
    
    The rebalancing interval time is assumed to be:
    \begin{equation}
        \hat{i}_{reb} = \frac{\hat{E}_{b,nom}}{\hat{P}_{b,nom}} \cdot 1.5
    \end{equation}
    
    This would be about the time needed to recharge the entire battery at rated charging power after the mixing, going from an oxidation state of 3.5+/3.5+ (can be imagined as $SoC$ = -50\%) to an oxidation state of 2+/5+ ($SoC$ = 100\%).
    
    The servicing maintenance time is considered negligible.

    \subsubsection{Model equations' set}
    
    The optimization problem can be summarized with the following equations set:
    objective function (\ref{FO}) has to be maximized, while complying with power constraints for the grid (\ref{curtailment}), (\ref{pg1}), (\ref{pg2}), power constraints for the battery (\ref{v1})-(\ref{v8}), state of charge and state of energy constraints (\ref{vsoc})-(\ref{vsoend}), performance bounding constraints (\ref{pccons})-(\ref{pdcons}) and rebalancing constraints (\ref{reb_cons1})-(\ref{reb_cons2}).
    All the constraints are applied for every hour $i \in T$, and for every day $d \in D$, .
    
    \subsection{Case studies}
    
    The optimization problem is applied to two different case studies:    
    \begin{enumerate}
        \item The first case analyzes the energy wholesale arbitrage in a renewable power plant, which sells energy at a volatile market price; an energy storage system is installed in an existing power plant, to maximize the profit from the energy exchange, according to the variations in day-ahead energy market price.
        \item The second case adds an energy storage system in a small power plant of a renewable energy community, with domestic users, to maximize the saving associated with purchasing energy from the grid, and promote the self-consumption of locally produced renewable energy. It is a case of renewable energy load shift.
    \end{enumerate}
    
    For both cases, two different types of renewable energy sources (RES) are investigated: wind and solar photovoltaic.
    The plants are both assumed to be installed in southern Italy.
    
    The technical and economic parameters used for the renewable plants are:
    
    \begin{table}[htb]
        \centering
        \begin{tabular}{cccccc}
        \toprule
        \textbf{Case N.} & \textbf{Name} & \textbf{Type} & \textbf{RES Size {[}kW{]}} & \textbf{Discount rate} & \textbf{Investment time {[}years{]}}
            \\ \midrule
        \multirow{2}{*}{1} & \multirow{2}{*}{Arbitrage} & Wind & $\num{1e4}$\\ \cline{3-4} &  & PV & $\num{1e4}$ & 8\% & 20
            \\ \midrule
        \multirow{2}{*}{2} & \multirow{2}{*}{Domestic} & Wind & 80
            \\ \cline{3-4} &  & PV & 180 & 6\% & 20
            \\ \bottomrule
        \end{tabular}
        \caption{Renewable plant installation characteristics for each case study}
        \label{impianti_RES}
    \end{table}

    % Note: the information contained in this table on the discount rate for the case is unnecessary, as CRF is not calculated anymore in this new version of the article. It could be eliminated.
    
    The wind and PV capacity factors are considered constant for both case studies, for simplicity, and rated respectively at 40.3\% and 17.9\%.
    The PV plant hourly production data are taken from the online PVGIS database and calculator \cite{JRC} and refer to a southern Italian test case. The wind plant hourly production data are taken from a study on a wind farm coupled with battery energy storage \cite{Frate2018}. Both renewable power plants are assumed to have a balance of plant efficiency of 88\%.
    
    For the first case study, the wholesale arbitrage, both the wind and PV plants sizes are 10 MW each.
    In this case, there is no electricity demand to satisfy and the electricity is only sold to the grid when it is most convenient, and never purchased from the grid.
    The selling price of electricity, in this case, is the market price on the day-ahead market. The reference prices are relative to the Sicilian market zone data from 2019 \cite{GME}.
    
    For the second case study, the domestic case, wind, and PV plants are sized according to the needs of the energy community, and the electricity can be both sold and purchased from the grid, at different prices.
    The hourly electricity demand of the renewable energy community is represented by employing historical data from 2019 reported by the European Network of Transmission System Operators for Electricity (ENTSO-e) \cite{ENTSO-e} through the convenient data collection provided by the Open Power System Data (OPSD) initiative \cite{TimeSeries}.
    The renewable energy community is assumed to be composed of 24 domestic users residing in the same neighborhood, and it has a peak consumption power of 72 kW.
    The load profile of the renewable energy community is assumed to be similar to the power consumed by the Sicilian bidding zone for the Italian market in 2019 scaled down by a \num{45000} factor.
    Both the wind and PV plants are sized to produce annually 70\% of the total energy annually consumed by the domestic user.
    The PV plant rated power is 180 kW and the wind plant rated power is 80 kW. Both RES plants and VRFB sizes are rounded up to the nearest multiple of five.
    For the sake of simplicity, the maximum power that can be exchanged with the grid for both cases equals twice the power of the installed renewable system, thus there are no connection limitations to the grid.
    In the domestic case, the energy is sold to the grid at a fixed price of $\hat{c}_{s} = 50 \EUR{}/MWh$ \cite{CER} and purchased from the grid at the mean price of $\hat{c}_{p} = 230 \EUR{}/MWh$ (Eurostat, 2019) \cite{Eurostat}.
    Furthermore, for the energy community, some incentives are applied to remunerate the self-consumption of energy and to promote the installation of electrochemical storage. The incentives on the Italian territory amount to 118 €/MWh of remuneration for every unit of self-consumed energy; this incentive is added to the fiscal detraction of 50\% of the investment cost, to be disbursed during the first 10 years of plant life \cite{CER}. Other European countries have similar incentives promoting renewable plants and energy storage installation in residential buildings and micro-grids, such as remunerating energy self-consumption or guaranteeing an appropriate feed-in tariff for renewable energy production \cite{Batstorm5}.
    
    The battery system has the following properties:
    
    \begin{table}[hbt]
        \centering
        \begin{tabular}{ccccc}
        \toprule
        \textbf{Case N.} & \textbf{Name} & \textbf{RES Type} & \textbf{VRFB Power {[}kW{]}} & \textbf{VRFB E/P ratio {[}h{]}} 
            \\ \midrule
        \multirow{2}{*}{1} & \multirow{2}{*}{Arbitrage} & Wind & 2500 \\ \cline{3-4} &  & PV & 2500 & 4
            \\ \midrule
        \multirow{2}{*}{2} & \multirow{2}{*}{Domestic} & Wind & 20
            \\ \cline{3-4} &  & PV & 45 & 4
            \\ \bottomrule
        \end{tabular}
        \caption{Battery storage system characteristics for each case study}
        \label{batt_size}
    \end{table}
    
    The energy coming from the RES plant is considered to be in form of direct current, entering the battery at a constant voltage. Then, an inverter is added to the battery system, to convert the output direct current into an alternate current, which is sent to the domestic users or sold to the grid. The inverter efficiency, variable with the output power, is added to the battery efficiency.

    \subsection{Economic aspects} \label{economic}
    
    Some economic aspects have to be considered to compare the results obtained using different battery models.
    
    \subsubsection{Maintenance costs calculation}
    
    In this subsection, we will present the economic aspects needed to calculate the total maintenance costs of the battery and the annual revenue associated with the battery use, which equals to the optimized value of the objective function, using Rodby's cost model \cite{Rodby2020}.
    
    To evaluate the annual O\&M costs, we sum the cost of every rebalancing and servicing during a whole year.
    
    A rebalancing event consists of mixing the two electrolytes by opening a valve between the two circuits for a few minutes to reach an equilibrium in the oxidation states, corresponding to a value of 3.5+ in both tanks.
    With this event, all the energy previously contained inside the battery is wasted, and the battery has to be recharged completely before the next use.
    The energy spent to charge the electrolytes to the maximum state of charge is made up of two contributions:
    (a) the energy spent to bring the oxidation state values from the homogeneous state of 3.5+/3.5+, to $SoC = 0\%$, corresponding to an oxidation state of 3+(anolyte)/4+(catholyte);
    (b) the energy spent to recharge the battery from $SoC$ = 0\% to the maximum feasible state of charge ($\hat{SoC}_{max}$).
    Using constraints (\ref{reb_cons1}) and (\ref{reb_cons2}), the energy spent to recharge the battery from the initial daily $SoC$ to $\hat{SoC}_{max}$ is already taken into account in the model, which also considers the time required for rebalancing ($\hat{i}_{reb}$), during which the battery can only be recharged, and not discharged.
    Additionally, we also have to account for the energy required to shift the oxidation states from 3.5+/3.5+ to 3+/4+, plus the energy required to shift from 0\% $SoC$ to the initial daily $SoC$ (from 3+/4+ to 2.7+/4.3+).
    This additional O\&M cost, not included in the optimization model, but included in the investment analysis, is:
    
    \begin{equation}
        C^{r}_{O\&M,a} = \sum_{d_{reb}=1} ^{N_R} C^{r}_{O\&M}(d) \Bigr|_{\substack{d = d_{reb}}} = \sum_{d_{reb}=1} ^{N_R} \hat{c}_{p}(1,d_{reb}) \cdot \frac{(E_m(d_{reb})\cdot 0.5 + \hat{SoE}_0)}{\bar{\eta}_c}
    \end{equation}
    Where $N_R$ is the annual number of rebalancing events, $\hat{c}_{p}$ is the cost of purchased energy from the grid, $\hat{SoE}_0$ is the value of SoE at the beginning of each day, and $\bar{\eta}_c$ is the mean efficiency of the battery during the charging operation for the rebalancing.
    We assume that the rebalancing is always performed at the beginning of the day, which is why the cost of purchasing energy from the grid refers to the first hour of the day.
    $\bar{\eta}_c$ is calculated assuming that the charging occurs at rated charging power and an average $SoC$ of 20\%.
    
    For chemical servicing maintenance, a certain amount of oxalic acid is added to the catholyte to reduce the vanadium pentoxide (V). The cost of this maintenance event is shown in the equations below, as elaborated by Rodby et al. \cite{Rodby2020}:
    
    \begin{equation}
        C^s_{O\&M} = E_b \cdot \left( O + \frac{n_{ox}\cdot MW_{ox} \cdot C_{ox}}{w_{ox}} \right)
    \end{equation}
    \begin{equation}
        n_{ox} = \left(U \cdot F \cdot \frac{\ 1 mol \: e-}{1 mol \: V}\right)^{-1} \cdot \frac{1 mol \: ox}{1 mol \: V} \cdot 3.6 \cdot 10^6 \frac{Ws}{kWh}
    \end{equation}
    
    Where $O$ is the labor cost, expressed in \$/kWh, $n_{ox}$ is the number of oxalic acid moles introduced for every kWh of installed energy, $MW_{ox}$ is the molecular weight of oxalic acid, $C_{ox}$ is its specific cost and $w_{ox}$ is its weight purity \cite{Rodby2020}.
    
    \begin{table}[hbt]
        \centering
        \begin{tabular}{cccc}
            \toprule
            Symbol & Description & Value & Unit \\
            \midrule
            O & Operational cost rate & 1 & $\$/kWh$ \\
            U & Open circuit voltage & 1.4 & $V$ \\
            F & Faraday's constant & \num{96485.33} & $C/mol e^-$\\
            $MW_{ox}$ & Molecular weight of oxalic acid & 90.03 & $g /mol$\\
            $C_{ox}$ & Cost of oxalic acid & 1.10 & $\$/ kg$\\
            $w_{ox}$ & Oxalic acid purity & 0.996 & - \\
            \bottomrule
        \end{tabular}
        \caption{Servicing related parameters name, description, unit, and value \cite{Rodby2020}}
        \label{Grandezze_OM}
    \end{table}
    
    The cost for every servicing event amounts to $3.65 \hspace{1mm} \$$ for every kWh of installed energy (corresponding to $3.0165 \EUR{}/kWh$ if we consider a 1.21 value for the dollar to euro exchange rate).
    The servicing cost model is assumed to be valid for the case where $cap_{lim}$ has a value of 80-85\%.
    
    The total O\&M annual costs are:
    
    \begin{equation}
        C_{O\&M,a} = C^r_{O\&M,a} + C^s_{O\&M,a} = C^r_{O\&M,a} + C^s_{O\&M}\cdot N_S \hspace{1cm} [\EUR{}/year]
    \end{equation}
    
    Where $N_S$ is the number of servicing events in a year, calculated by the optimization model.
    
    \subsubsection{Net annual revenue calculation}
    
    In this subsection, we calculate the net annual revenue associated with battery use.
    The net revenue is calculated as a difference between the revenue (optimized value of the objective function) associated with the storage system use and the revenue that we would have in the same case study, but without the VRFB:
    
    \begin{equation}
        \Delta R_{a} = Rev_a - Rev_{a,wob} + \Delta I_a \hspace{1cm} [\EUR{}/year]
    \end{equation}
    
    $Rev_a$ is the annual revenue associated with electricity selling and purchasing from the grid when the use of VRFB storage is optimized.
    
    \begin{equation}
        Rev_a = \sum_{d=1}^{365} Rev(d)
    \end{equation}
    
    $Rev_{a,wob}$ is the annual revenue associated with electricity selling and purchasing from the grid, in a hypothetical scenario where the RES plant is installed without the battery system and the energy is directly sold to the grid, following the production curve.
    
    \begin{equation}
        Rev_{a,wob} = \sum_{d=1}^{365} \tau \cdot \sum_{i=1}^{24}
        \begin{cases}
            \hat{c}_s(i,d) \cdot (\hat{P}_{res}(i,d) - \hat{P}_{dem}(i,d)) & \text{if} \hspace{1 mm} \hat{P}_{res}(i,d) \geq \hat{P}_{dem}(i,d) \\
            \hat{c}_p(i,d) \cdot (\hat{P}_{dem}(i,d) - \hat{P}_{res}(i,d)) & \text{if} \hspace{1 mm} \hat{P}_{res}(i,d) < \hat{P}_{dem}(i,d) \\
        \end{cases} 
    \end{equation}
    
    $\Delta I_a$ is a term that accounts for the incentive tariff applied to the self-consumption of energy by the Italian regulation.
    $\Delta I_a$ is zero for the first case study, as there are no incentives for the energy arbitrage purpose, but is calculated as follows for the second case study:
    
    \begin{equation}
        \Delta I_a = 
        \begin{cases}
            I_a - I_{a,wob} + C_{a,det} & \text{if} \hspace{1 mm} a \in [1,10] \\
            I_a - I_{a,wob}             & \text{if} \hspace{1 mm} a \in [11,20]\\
        \end{cases} 
    \end{equation}
    
    $\Delta I_a$ changes during the investment lifetime, considering that part of the incentive is made partly of the remuneration of self-consumed energy and partly of the fiscal detraction of 50\% of the investment cost for the battery system ($C_{det}$).
    
    The incentives terms for the remuneration of self-consumed energy are evaluated as follows, in the case with ($I_a$) and without ($I_{a,wob}$) the battery:
    
    \begin{equation}
        I_a = C_{inc} \cdot \tau \cdot \sum_{d=1}^{365} \sum_{i=1}^{24} P^{sc}(i,d) \hspace{0.5 cm} [\EUR{}/year]
    \end{equation}
    \begin{equation}
        I_{a,wob} = C_{inc} \cdot \tau \cdot \sum_{d=1}^{365} \sum_{i=1}^{24} P^{sc}_{wob}(i,d) \hspace{0.5 cm} [\EUR{}/year]
    \end{equation}
    
    $C_{inc} = 118 \EUR{}/MWh$ is the remuneration incentive yielded for every unit of energy directly consumed; while $P^{sc}(i,d)$ is the amount of power that is directly self-consumed during hour \textit{i} and day \textit{d}, with the battery installed, and $P^{sc}_{wob}(i,d)$ is the amount of power that is directly self-consumed, without the battery installed.
    The self-consumed powers are calculated as follows:
    
    \begin{equation}
        P^{sc}(i,d) =
        \begin{cases}
            \hat{P}_{dem}(i,d) & \text{if} \hspace{1 mm} k_{id} (i,d) = 0\\
            \hat{P}_{dem}(i,d) - P_{g,p}(i,d) & \text{if} \hspace{1 mm} k_{id} (i,d) = 1 \\
        \end{cases}
    \end{equation}
        
    \begin{equation}
        P^{sc}_{wob}(i,d) = min(\hat{P}_{res}(i,d),\hat{P}_{dem})
    \end{equation}
    
    The incentive of fiscal detraction of 50\% of the investment cost for the battery system ($C_{a,det}$), disbursed during the first 10 years of the investment, is calculated as follows:
    
    \begin{equation}
        C_{a,det} = 50\% \cdot \frac{C_{b}}{10}
    \end{equation}
    
    Where $C_{b}$ is the capital cost of the battery system and equals the sum of the total power-related costs and the total energy-related costs:
    
    \begin{equation}
        C_b = \hat{P}_{b,nom} \cdot \hat{C}_P + \hat{E}_{b,nom} \cdot \hat{C}_E \hspace{1cm} [\EUR{}]
    \end{equation}
    
    The specific battery costs in terms of power and energy installed are $\hat{C}_P = 1080 \EUR{}/kW$ and $\hat{C}_E = 385 \EUR{}/kWh$  \cite{Minke2017}.
    
    Note that for the second case study, the mean value of the annual net revenue $\Delta R$ during the whole battery operational life (20 years) is the mean value between the first 10 years of the investment (with the fiscal detraction incentives) and the last 10 years of the investment (without the fiscal detraction incentive):
    
    \begin{equation}
        \Delta R = \sum_{a=1}^{20} \frac{\Delta R_a}{20}
    \end{equation}
    
    \subsection{Simple models}
    The detailed optimization model is compared, in terms of energetic and economic results, to two simpler models, based on the simplification commonly used in the literature to evaluate the errors given by an incorrect description of the system.
    
    \begin{enumerate}
        \item the first simple model uses efficiencies as a function of power and $SoC$, but does not take into account the degradation effects and maintenance;
        \item the second simple model uses fixed efficiencies and also does not take into account the degradation effects and maintenance.
    \end{enumerate}
    
    In both simple models, the battery does not experience any degradation or maintenance.
    Furthermore, to fairly compare simple and detailed models, the simple model with fixed efficiencies uses the mean efficiency values calculated using the detailed model: $\hat{\eta}^*_{c}$ and $\hat{\eta}^*_{d}$.
    
    \bigskip
    Table \ref{parameters_values} contains names, descriptions, units of measurement and values of all the parameters necessary to recreate and solve the optimization model.

    \begin{table}[htb]
    \centering
        \begin{tabular}{llll}
        \toprule
        \textbf{Name}     & \textbf{Description}                    & \textbf{Unit}  & \textbf{Value}   \\
        \midrule
        $\bar{\eta}_c$    & Mean efficiency during rebalancing      & {[}\%{]}       & 79.7             \\
        $\bar{\eta}^*_c$  & Charging efficiency for simple models   & {[}\%{]}       & 75.9             \\
        $\bar{\eta}^*_d$  & Discharging efficiency for simple models & {[}\%{]}       & 73.5             \\
        $\hat{B}M$                     & Big value to deactivate the constraints         & {[}-{]}                      & 1.5 $\cdot \hat{P}_{b,nom}$ \\
        $\hat{C}_E$       & Energy related battery price            & {[}€/kWh{]}    & 385              \\
        $\hat{C}_P$       & Power related battery price             & {[}€/kW{]}     & 1080             \\
        \multirow{2}{*}{$\hat{c}_{p}$} & \multirow{2}{*}{Price of purchased electricity}      & \multirow{2}{*}{{[}€/kWh{]}} & 0 for Case 1         \\
                          &                                         &                & 0.230 for Case 2 \\
        \multirow{2}{*}{$\hat{c}_{s}$} & \multirow{2}{*}{Price of sold electricity} & \multirow{2}{*}{{[}€/kWh{]}} & $f(i,d)$ for Case 1  \cite{GME}    \\
                          &                                         &                & 0.050 for Case 2 \\
        $\hat{cap}_{lim}$ & Lower capacity limit                    & {[}\%{]}       & 80               \\
        $n_{int}$         & Number of linearization intervals       & {[}-{]}        & 5                \\
        $\hat{P}_{g,max}$              & Maximum power exchanged with the grid           & {[}-{]}                      & 2 $\cdot \hat{P}_{nom,res}$ \\
        $\hat{R}_{fade}$  & Capacity fade rate                      & {[}\%/cycle{]} & 0.442            \\
        $\hat{r}_{ED}$    & Electrolyte decay rate                  & {[}\%/cycle{]} & 0.055            \\
        $\hat{SoC}_0$     & Initial state of charge                 & {[}-{]}        & 0.3              \\
        $\hat{SoC}_{max}$ & Maximum state of charge                 & {[}-{]}        & 0.9              \\
        $\hat{SoC}_{min}$ & Minimum state of charge                 & {[}-{]}        & 0.1              \\
        \bottomrule
        \end{tabular}
        \caption{Name, description, unit of measurement, and values of the parameters used in the optimization model}
        \label{parameters_values}
    \end{table}
    
    \section{Results and discussion} \label{results}
    
    This section presents and discusses the most relevant technical and economical results associated with the optimal management strategy of the VRFB system described above for the presented case studies.
    
    The problem is solved using a computer with the following specifics: CPU Intel(R) Core(TM) i7-8550u with a 1.80 GHz clock speed, with 4 physical cores and 8 logical processors.
    The daily optimization problem is solved, inside the Yalmip \cite{yalmip} environment for Matlab, using the Gurobi \cite{gurobi} solver, within an average time of 8.22 s for the domestic case and 4.93 s for the arbitrage case.
    
    \subsection{Capacity decay prediction with detailed model}
    
    In this section, we compare the results obtained using the optimization algorithm developed in the paper, with the results obtained with the simplifying hypothesis of constant daily battery cycling, which is a common hypothesis in techno-economic studies.
    
\begin{figure}[hbt]
        \centering
        \subfigure[]{%
        \includegraphics[scale = 0.58]{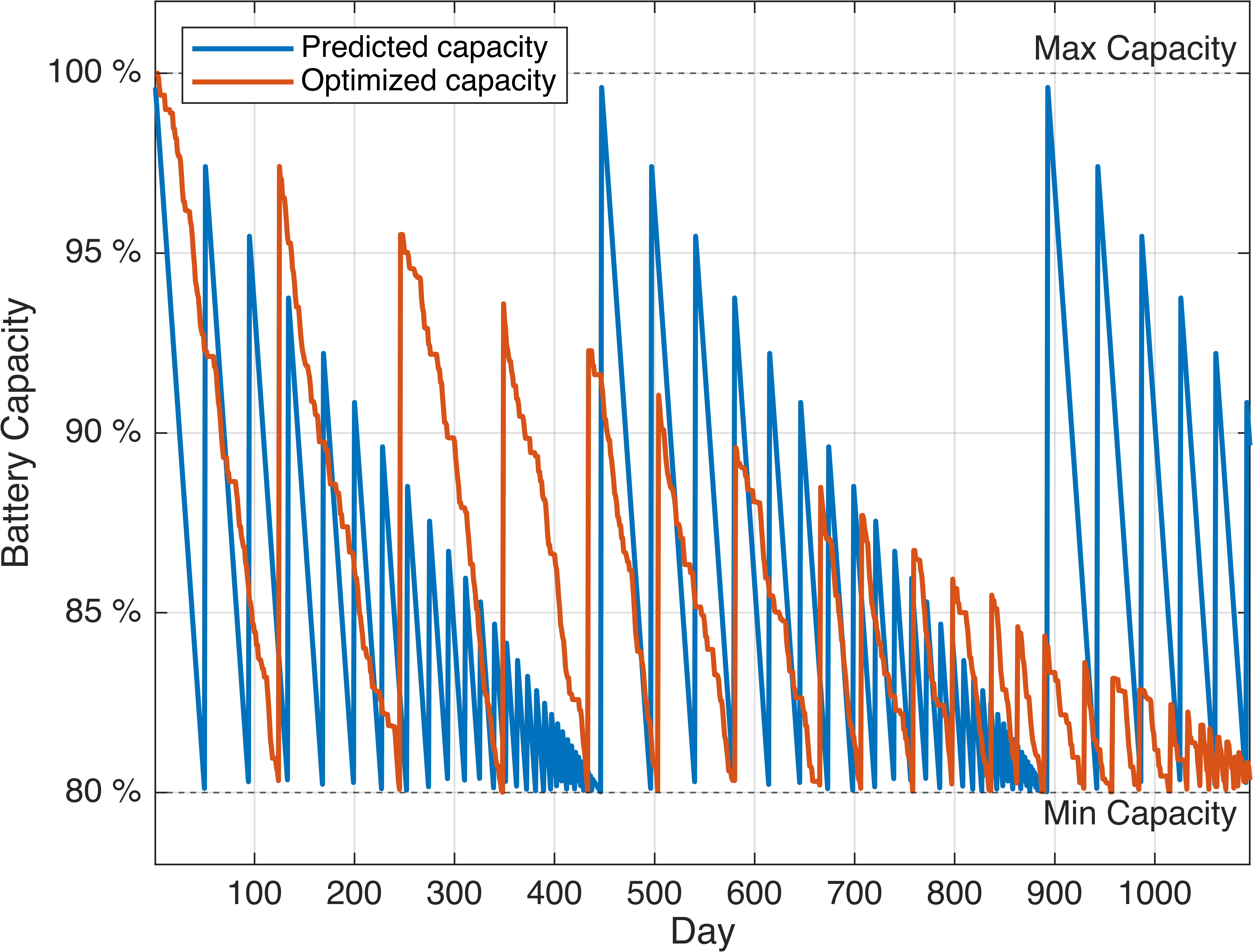}
        }
        \quad
        \subfigure[]{%
        \includegraphics[scale = 0.58]{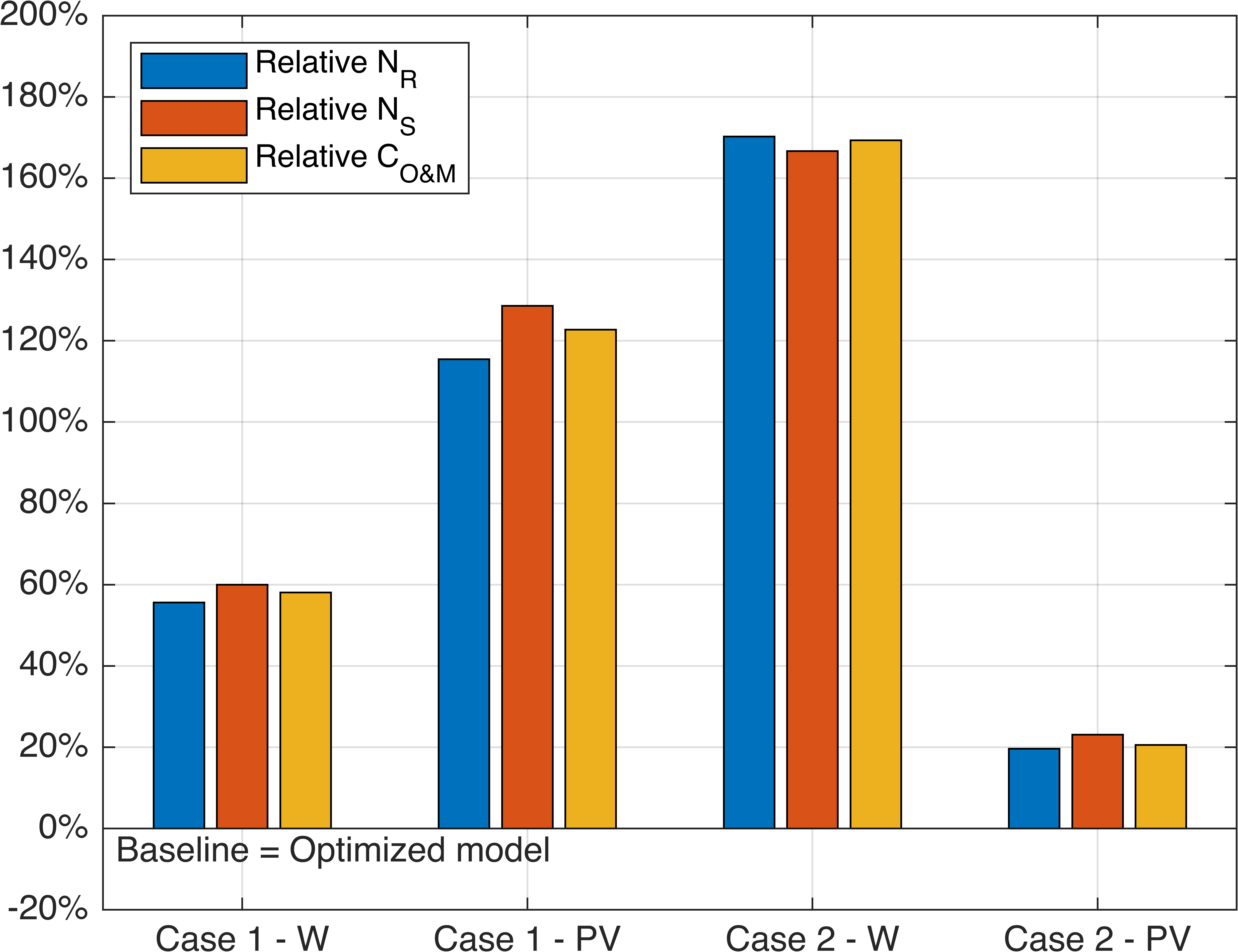}
        }
        \caption{Optimized vs. predicted results: (a) comparison of capacity decay and maintenance during three years of operation, case 2-W; (b) relative difference between predicted and optimized (baseline value) number of maintenance events and costs, during battery operational lifetime, for arbitrage case (1) and domestic case (2), with wind plant (W) and photovoltaic plant (PV)}
        \label{cap_decay_comparison}
    \end{figure}

    Figure \ref{cap_decay_comparison}(a) shows the fraction of accessible capacity (relative to the rated capacity value) for three years of continuous operation, for one representative case (2-W: domestic case study, with a wind plant installed).
    Each year is assumed to have the same features as the others (same energy price, production, and consumption profiles).
    The figure shows the comparison between the optimized fraction of accessible capacity, calculated using the program presented in this paper, and an ideal fraction of capacity, predicted assuming that the battery completes a constant number of daily charge/discharge cycles.
    The simplifying hypothesis for a constant number of charge/discharge cycles is the following:
    the battery completes one charge/discharge cycle per day, within its state of charge limits, or at a certain DoD (depth of discharge) range.
    The DoD range is defined as follows:
    
    \begin{equation}
        DoD = \hat{SoC}_{max} - \hat{SoC}_{min}
    \end{equation}
    
    The equivalent number of charge/discharge cycles per day, in the predicted model, is:
    
    \begin{equation}
        N_{cyc,d}^* = 1 \cdot DoD
    \end{equation}
    
   It results here that $ N_{cyc,d}^* = DoD = 0.8 $.
    
    Assuming to evaluate the whole operational lifetime of the battery (20 years), the number of maintenance events, and total cost, are represented in Table \ref{comparison_maintenance} with their actual values, and in Figure \ref{cap_decay_comparison}(b) with their relative values.
    $N_R$ is the total number of rebalancing events, $N_S$ is the total number of servicing events, $C_{O\&M}$ is the total operation and maintenance cost during the lifetime. To evaluate $C_{O\&M}$ for the predicted case, a variation of the rebalancing cost, presented in subsection \ref{economic}, has been employed. The rebalancing cost formula is built considering the annual mean values of the following variables subject to the summation: $\hat{c}_p(i,d)$ and $E_m(d)$. This allows obtaining the mean cost for a single rebalancing operation, which can be multiplied by the total number of predicted rebalancing operations to obtain the total cost and to correctly compare the two results, calculated with the same formula.
    
    \begin{table}[htb]
        \centering
        \begin{tabular}{cccc}
            \toprule
            \textbf{Case} &
              \textbf{\begin{tabular}[c]{@{}c@{}}Parameter\\ name\end{tabular}} &
              \textbf{\begin{tabular}[c]{@{}c@{}}Optimized\\ value\end{tabular}} &
              \textbf{\begin{tabular}[c]{@{}c@{}}Predicted\\ value\end{tabular}} \\
            \midrule
            \multirow{3}{*}{1-W}  & $N_R$      & 402                & 627               \\
                                  & $N_S$      & 10                 & 16                \\
                                  & $C_{O\&M}$ & 5.35 $\cdot 10^5 $ & 8.47 $\cdot 10^5$ \\
            \midrule
            \multirow{3}{*}{1-PV} & $N_R$      & 291                & 627               \\
                                  & $N_S$      & 7                  & 16                \\
                                  & $C_{O\&M}$ & 3.80$\cdot 10^5$   & 8.48 $\cdot 10^5$ \\
            \midrule
            \multirow{3}{*}{2-W}  & $N_R$      & 232                & 627               \\
                                  & $N_S$      & 6                  & 16                \\
                                  & $C_{O\&M}$ & 5.28 $\cdot 10^3$  & 1.42 $\cdot 10^4$ \\
            \midrule
            \multirow{3}{*}{2-PV} & $N_R$      & 524                & 627               \\
                                  & $N_S$      & 13                 & 16                \\
                                  & $C_{O\&M}$ & 2.63 $\cdot 10^4$  & 3.17 $\cdot 10^4$ \\
            \bottomrule
            \end{tabular}
        \caption{Optimized vs. predicted number of maintenance events [cycles] and costs [€] during battery operational lifetime (20 years), for arbitrage case (1) and domestic case (2), with wind plant (W) and photovoltaic plant (PV)}
        \label{comparison_maintenance}
    \end{table}
    
    Figure \ref{cap_decay_comparison}(b) shows the relative difference between the predicted and optimized values for the examined variables (number of rebalancing events, servicing events, operation and maintenance costs) and calculated during the battery operational lifetime, for different cases.
    Assuming a constant number of charging/discharging cycles per day throughout the whole year can result in a substantial error in the estimation of the capacity fade and the accessible capacity of the battery if compared to the results obtained for the optimized model.
    This error in the capacity fade prediction could result in errors in maintenance events scheduling. This will then lead to errors in the evaluation of the total $O\&M$ costs.
    The figure highlights that the greatest problem with the prediction of capacity fade is the overestimation of the battery cycling, which yields an overestimation of the total expenses of the battery.
    This occurs especially in cases where the optimized use of the battery is modest and the battery performs a low number of charging/discharging cycles per day.
    The aforementioned cases are the arbitrage case (Case 1), with PV plant, where the mean annual number of cycles per day is $N_{cyc,d} = 0.38$, and the domestic case (Case 2), with wind plant, where the mean annual number of cycles per day is $N_{cyc,d} = 0.31$.
    For these cases (Case 1 - PV and Case 2 - W), the number of rebalancing events predicted with the constant number of cycling hypothesis, are respectively, 115\% and 170\% higher than the number of rebalancing calculated with the optimized model; the predicted number of servicing events are 129\% and 167\% higher, and the total O\&M costs are 123\% and 169\% higher than in the optimized model.
    For the other cases (Case 1 - W and Case - PV), the differences between the predicted and optimized model are less relevant, up to 60\% for the arbitrage case (Case 1) with wind plant. These cases have higher mean annual number of cycles per day, which are: $N_{cyc,d} = 0.49$ for Case 1 - W, and $N_{cyc,d} = 0.67$ for Case 2 - PV.
    The results concerning the optimal scheduling of maintenance events are highly dependent on the selected test case study.
    
    \subsection{Global impact of detailed model versus simple models}
    The following section shows the comparison between the techno-economic results obtained with the detailed optimization model, which uses variable efficiencies and capacity fade with restoration, and the results given by two basic models, which do not take into account the battery efficiency curves and degradation. It appears that different models lead to different use strategies for the battery and different economical revenue.
    The results are shown for different case studies, with different renewable plants.
    
\begin{figure}[hbt]
        \centering
        \subfigure[]{%
        \includegraphics[scale = 0.6]{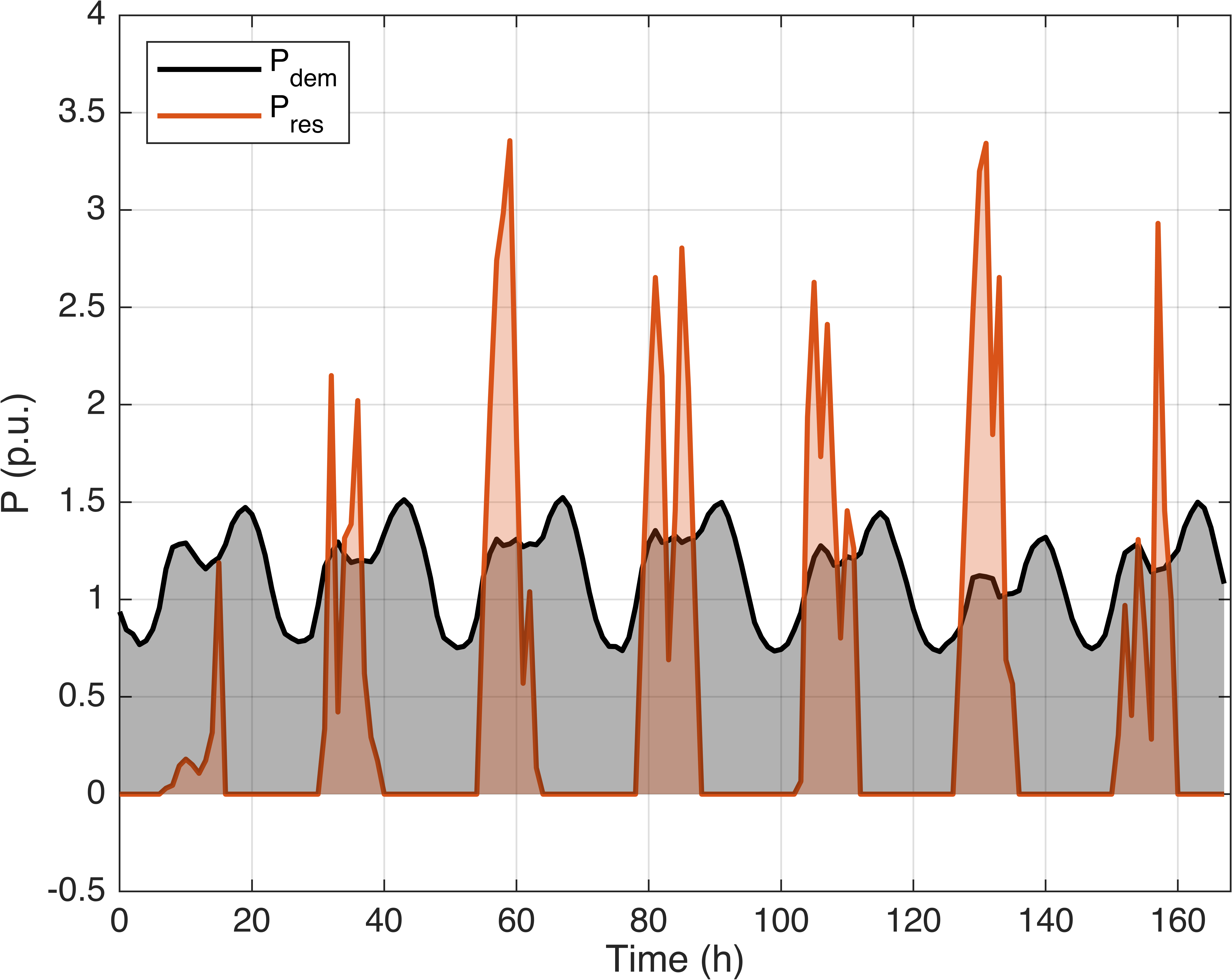}
        }
        \quad
        \subfigure[]{%
        \includegraphics[scale = 0.6]{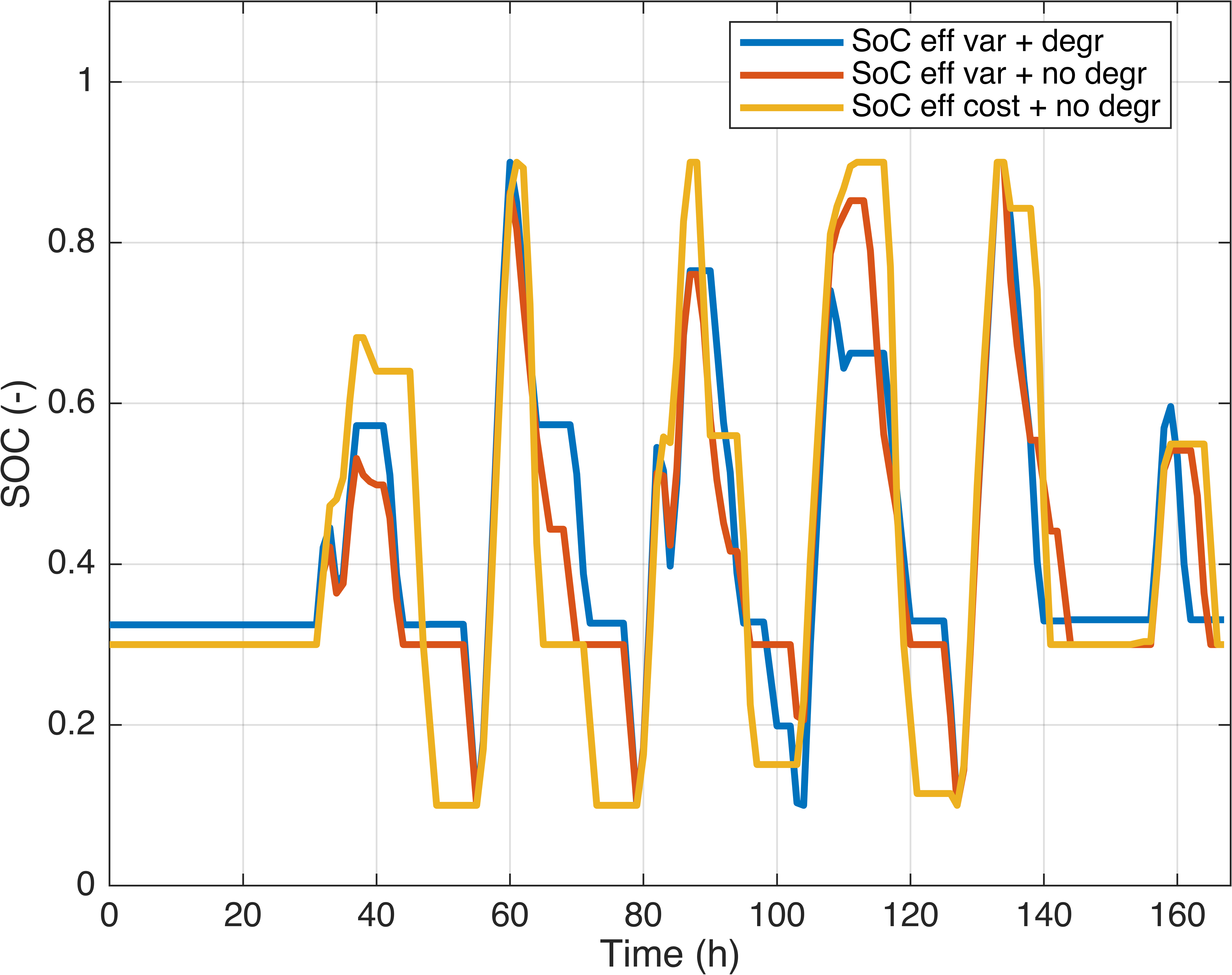}
        }
        \caption{Winter week analysis: dimensionless produced and demanded power (a) and $SoC$ comparison for three models (b), for the domestic case study, with PV plant}
        \label{midseason_comp}
    \end{figure}

    Figure \ref{midseason_comp} shows the input data and the results for the domestic case study, with a photovoltaic plant, during a winter week ($4^{th}$ week of the year).
    Figure \ref{midseason_comp}(a) plots the dimensionless values (divided by the rated battery power) of demanded power from the renewable energy community versus the renewable production.
    Figure \ref{midseason_comp}(b) shows the battery management results in the form of $SoC$, for the three different models.
    Different models give slight differences in output due to the different efficiencies and capacity values, but overall the management is similar because it is induced by demanded and produced renewable disparities during the day, as can be seen in Figure \ref{midseason_comp}(a).
    
\begin{figure}[hbt]
        \centering
        \subfigure[]{%
        \includegraphics[scale = 0.6]{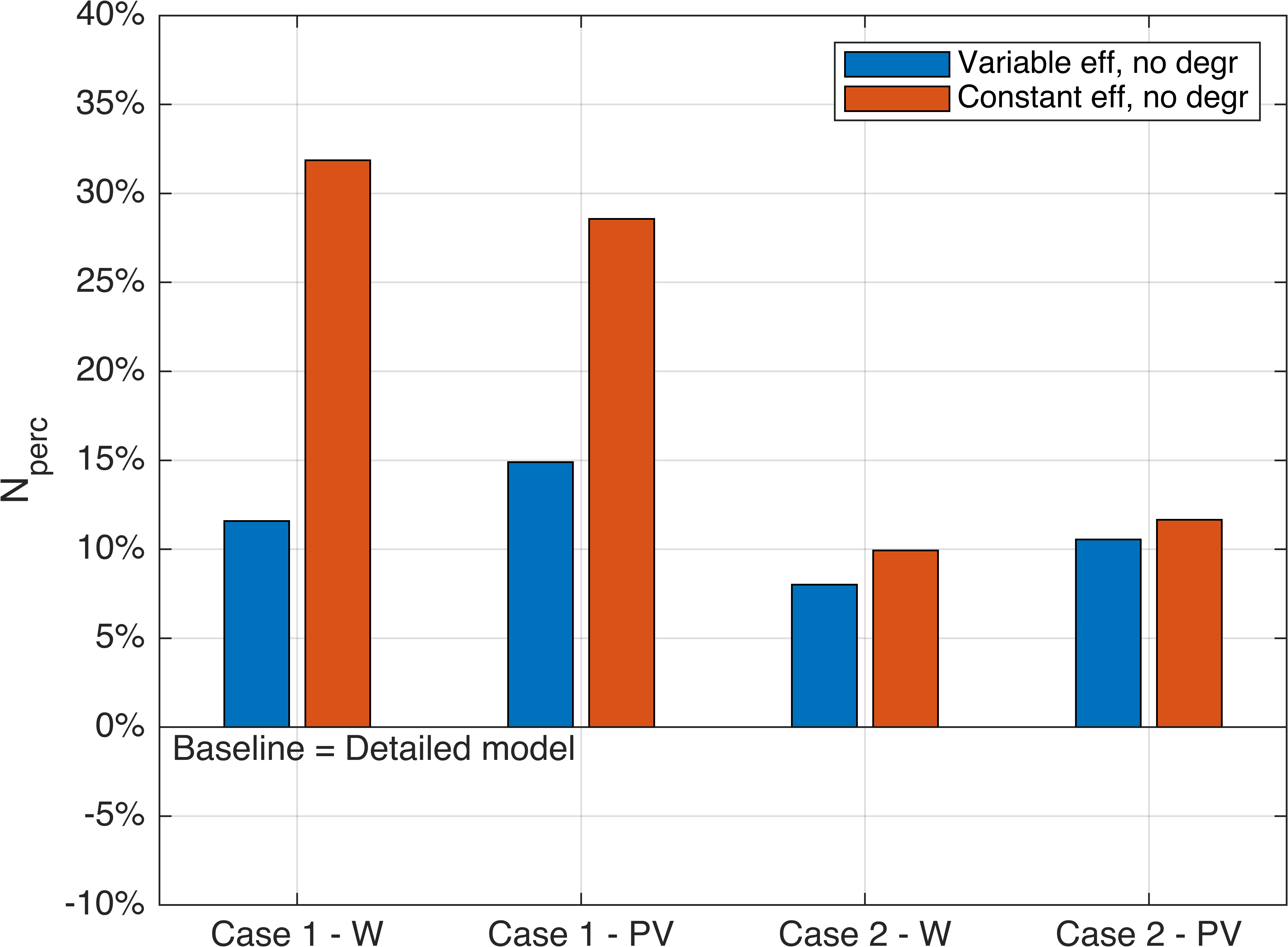}
        }
        \quad
        \subfigure[]{%
        \includegraphics[scale = 0.6]{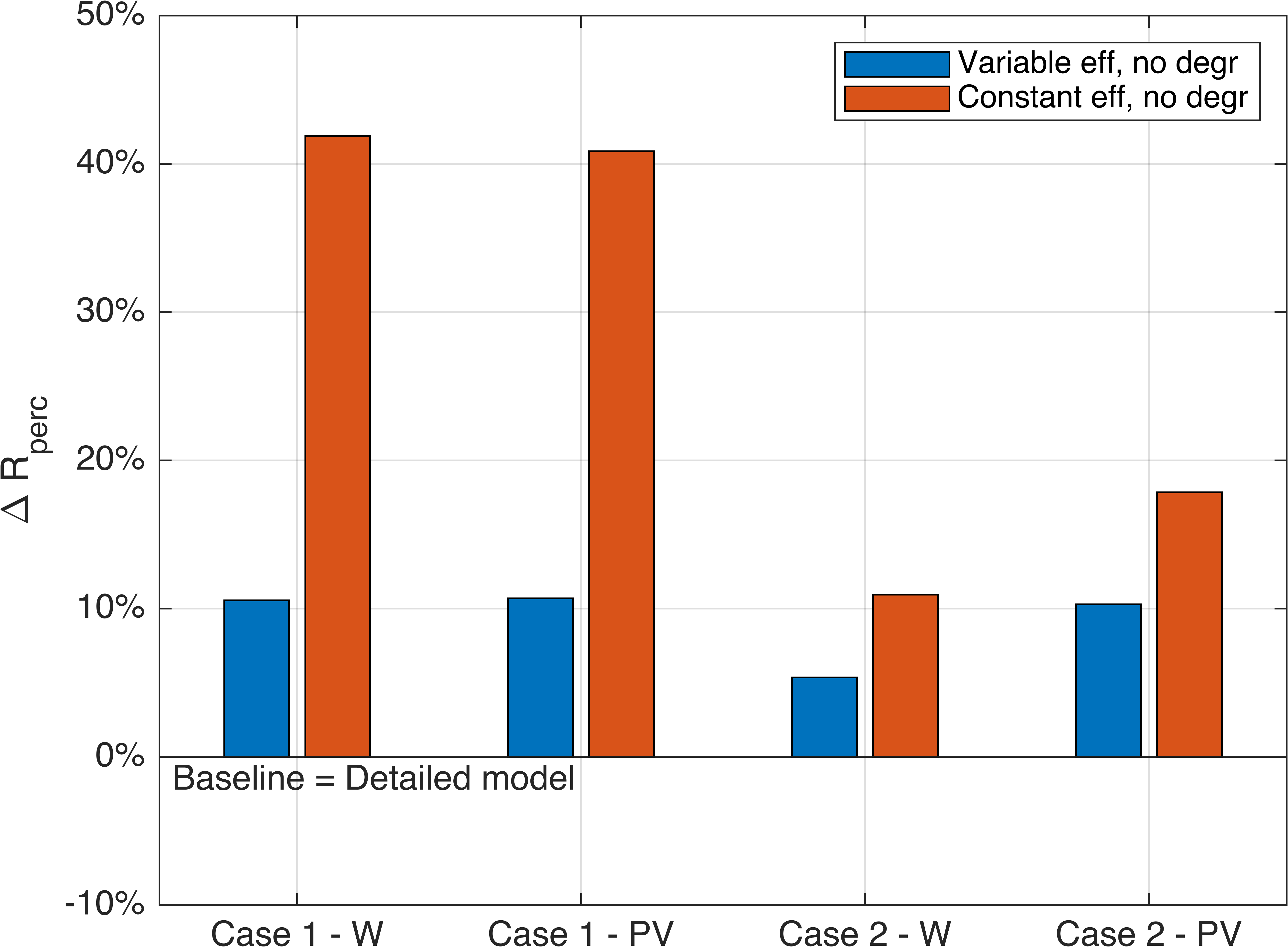}
        }
        \caption{Annual relative differences between simple and detailed (baseline value) models for arbitrage case (1) and domestic case (2), with wind plant (W) and photovoltaic plant (PV): (a) number of cycles; (b) net revenue}
        \label{confmodel}
    \end{figure}

    As far as the annual number of charge-discharge cycles of the battery, those are compared, for all the case studies, in Figure \ref{confmodel}(a).
    On the y-axis, there is the percentage deviation of the annual number of cycles:
    
    \begin{equation}
        N_{perc} = \frac{N^0 - N}{N}\cdot 100 \hspace{0.5 cm} [\%]
    \end{equation}
    
    Where $N$ is the annual number of cycles calculated with the detailed model; while $N^0$ is the annual number of cycles calculated with any simple model.

    For the first case study, if the capacity fade of the battery is neglected, the number of ideal cycles of the battery is 11\% to 15\% higher than with the detailed model, depending on the type of renewable plant (W = wind, PV = photovoltaic), while if there are both constant efficiency and no degradation model, the number of ideal cycles increases by a total of 29\% to 32\%, when compared with the detailed model.
    For the second case study, if the capacity fade of the battery is neglected, the number of ideal cycles is 8\% to 11\% higher than with the detailed model, while if there are both constant efficiency and no degradation, the number of cycles increases by a total of 10\% to 12\%, when compared with the detailed model.
    The model that does not take into account the capacity fade overestimates the availability of the energy stored in the battery: this results in a higher number of ideal charge and discharge cycles.
    Furthermore, if the efficiencies are considered constant, there is no penalization due to different power and $SoC$, and this will also underestimate the energy losses, leading to more battery utilization.
    
    Additionally, simple models lead to an overestimation of the revenue, as can be seen in Figure \ref{confmodel}(b), which shows the percentage deviation of the annual net revenue (or saving) calculated with the simple models, in comparison to the detailed model, for the different case studies.
    On the y-axis, the percentage deviation of the revenue is shown:
    
    \begin{equation}
        \Delta R_{perc} = \frac{\Delta R^0 - \Delta R}{\Delta R}\cdot 100 \hspace{0.5 cm} [\%]
    \end{equation}
    
    Where $\Delta R$ is the mean value of the annual net revenue during the whole battery operational life, with the detailed model; while $\Delta R^0$ is the mean value with any simple model.
    
    For the first case study, if the capacity fade of the battery is neglected, the annual revenue of the battery is 11\% higher than with the detailed model for both renewable plants, while if there are both constant efficiency and no degradation model, the number of ideal cycles increases by a total of 41\% to 42\%, when compared with the detailed model.
    For the second case study, if the capacity fade of the battery is neglected, the number of ideal cycles of the battery is 5\% to 10\% higher than with the detailed model, while if there are both constant efficiency and no degradation model, the number of ideal cycles increases by a total of 11\% to 18\%, when compared with the detailed model.
    
    This is a direct consequence of the overestimation of the stored energy and number of cycles, driven by the poor calculation of energy losses with simple models.
    
    $N_{perc}$ and $\Delta R_{perc}$ overestimation with simple models is higher for the first case study than for the second case study. This happens because the values of $N_0$ and $\Delta R_0$, for the first case study, are relatively small, thus a small deviation in the model can lead to big relative errors.

    \section{Conclusions} \label{conc}
    
    % 1.
    Vanadium Redox Flow Batteries, with their high chemical stability, long operational life, and the possibility to restore the capacity fade with periodic maintenance without capital expenditure, are a promising technology for renewable energy storage and have possible applications in both residential and arbitrage cases.
    The scope of the paper was to define a detailed operation optimization model for a Vanadium Redox Flow Battery system, to better understand the effects of accurate modeling on the optimal battery management strategies, and ultimately to understand the limitations of the technology when paired up with renewable power sources in different test case studies.
    
    % 2.
    The detailed optimization model for VRFB scheduling was developed using variable efficiencies, functions of charging/discharging powers and state of charge, also accounting for capacity fade due to different electrolyte degradation mechanisms and adding different maintenance strategies to restore the battery capacity after periodic use. The capacity restoration model takes into account both costs and operational time of the maintenance events (rebalancing and chemical servicing), needed to reverse the effects of the two main causes of capacity fade: crossover and oxidative imbalance.
    The overall problem was solved daily, with an hourly resolution, as a MILP, while the nonlinearity of the efficiency functions was solved with the use of a two-variables convexification algorithm, yielding a faster algorithm, solved with low computational time.
    
    The model was applied to two different case studies, within the energy arbitrage and the domestic framework, using as a reference a south Italian bidding zone and its data for energy prices, renewable energy production curves, and typical regional household energy load. The renewable energy plants in the evaluated scenarios are either wind or photovoltaic.
    The optimization consisted of maximizing the total daily revenue from trading electricity with the grid.
    
    % 3.
    After developing this detailed model, the results were compared to those obtained with simpler battery models, which estimate an ideal use of the battery, by setting a fixed number of daily charge/discharge cycles, or which assume constant battery efficiency and do not consider capacity fade.
    
    Correctly optimizing the number of charging/discharging cycles in different scenarios has been proved to be challenging and could result in overestimation of the optimal number of maintenance events and costs, up to 170\%, if the charging/discharging cycles are assumed to be constant at one cycle per day at rated DOD range.
    
    Using simpler optimization models, which do not take into account the effect of variable efficiencies and capacity fade, leads to the overestimation of the optimal number of cycles of the battery by up to 15\%, and of the revenue by up to 11\%, depending on the scenario, if they do not take into account the degradation model of the battery, and respectively up to 32\% and 42\%, if they also assume constant efficiency for the battery when compared to the detailed optimization model presented in this work.

    \section*{CRediT authorship contribution statement}
    
    \textbf{Diana Cremoncini}
    Conceptualization, Methodology, Software, Data Curation, Visualization, Formal Analysis, Writing - Original Draft.
    \textbf{Guido F. Frate}
    Methodology, Software, Data Curation, Formal Analysis, Supervision, Writing - Review \& Editing.
    \textbf{Aldo Bischi}
    Conceptualization, Methodology, Supervision, Funding acquisition, Writing - Review \& Editing.
    \textbf{Lorenzo Ferrari} Supervision,  Project administration, Funding acquisition, Writing - Review \& Editing.
    
    \section*{Declaration of competing interest}
    The authors declare that they have no known competing financial interests or personal relationships that could have appeared to influence the work reported in this paper.

    \section*{Acknowledgements}
    The authors thank Dr. Andrea Baccioli for the supervising process and professor Antonio Bertei for many helpful discussions on the electrochemical processes concerning flow batteries.
    
    \section*{Funding}
    This research has received funding from the European Union’s Horizon 2020 research and innovation program under Grant Agreement No 875565 (CompBat: Computer-aided design for next generation flow batteries. Project website: https://compbat.eu/).
    
    This research is supported by the Ministry of University and Research (MUR) as part of the PON 2014-2020 “Research and Innovation" resources – Green/Innovation Action – DM MUR 1061/2022.
    
    Dr. Guido Francesco Frate acknowledges the financial contribution received from the Italian Operative National Plan (\textit{Piano Operativo Nazionale}, PON) in the framework of the project “\textit{Ricerca e Innovazione}” 2014-2020 (PON R\&I) \textit{Azione IV.4 “Dottorati e contratti di ricerca su tematiche dell'innovazione” (Azione IV.6 “Contratti di ricerca su tematiche Green”)}.
    
    \section*{Nomenclature}
    
    \begin{table}[H]
        \begin{tabular}{llc}
        \textbf{Symbols}    & \textbf{Description}              & \textbf{Unit} \\
        \midrule
        $\gamma$            & charging tangent plane coefficient &     -        \\
        $\delta$            & discharging tangent plane coefficient &  -        \\
        $\eta$              & efficiency                        & -             \\
        $\tau$              & time step length                  & h             \\
        BM                  & big value                         & -             \\
        C                   & capital cost of battery           & €/kW, €/kWh   \\
        c                   & price of electricity              & €/kWh         \\
        cap                 & capacity fraction                 & \%            \\
        d                   & d-th day                          & -             \\
        D                   & number of days                    & -             \\
        E                   & energy                            & kWh           \\
        f                   & generic factor                    & -             \\
        i                   & i-th time step                    & -             \\
        k                   & binary variable                   & -             \\
        n                   & numeric parameter                 & -             \\
        P                   & power                             & kW            \\
        R, r                & capacity fade rate                & \%/cycle      \\
        Rev                 & economic revenue                  & €/kWh         \\
        SoC                 & state of charge                   & -             \\
        SoE                 & state of energy                   & -             \\
        T                   & number of optimization time steps & -             \\
        \end{tabular}
    \end{table}

    % Table 2.B
    
    \begin{table}[H]
        \begin{tabular}{ll}
        \textbf{Subscript/superscript} & \textbf{Description}           \\
        \midrule
        0                              & initial state, reference state \\
        b                              & battery, battery state         \\
        c                              & charge                         \\
        cap                            & capacity                       \\
        curt                           & curtailment                    \\
        cyc                            & cycle number                   \\
        d                              & discharge, day                 \\
        dem                            & demanded                       \\
        ED                             & Electrolyte decay              \\
        fade                           & Total fade rate                \\
        g                              & grid                           \\
        id                             & grid state                     \\
        in                             & internal                       \\
        int                            & interval                       \\
        j                              & j-th tangent plane             \\
        lim                            & lower limit                    \\
        losses                         & losses due to inefficiencies   \\
        m                              & maximum stored                 \\
        max                            & maximum                        \\
        min                            & minimum                        \\
        nom                            & nominal                        \\
        onoff                          & battery state ON/OFF           \\
        p                              & purchase                       \\
        q                              & sampling point                 \\
        R, reb                         & rebalancing                    \\
        res                            & renewable energy source        \\
        S                              & servicing                      \\
        s                              & sell                          
        \end{tabular}
    \end{table}

\bibliographystyle{unsrt}
\bibliography{main}

\end{document}